\documentclass{article}
\usepackage{arxiv}
\usepackage{authblk}
\usepackage[utf8]{inputenc}
\usepackage{textcomp}
\usepackage{color}
\usepackage{graphicx}
\usepackage{float}
\usepackage{caption,subcaption}
\usepackage{amsmath}
\usepackage[version=4]{mhchem}
\usepackage{mathtools}
\usepackage{siunitx}
\usepackage{longtable,tabularx}
\usepackage{tikz}
\usepackage[many]{tcolorbox}
\usepackage{color}
\usepackage{epstopdf}
\usepackage[T1]{fontenc}    
\usepackage{hyperref}       
\usepackage{url}            
\usepackage{booktabs}       
\usepackage{amsfonts}       
\usepackage{nicefrac}       
\usepackage{microtype}      
\usepackage{lipsum}
\usepackage{amssymb}
\usepackage{algorithm}
\usepackage{tcolorbox} 
\graphicspath{ {./images/} }

\usepackage[numbers,sort&compress]{natbib}
\def\1{\boldsymbol{1}}

\newtheorem{problem}{Problem}

\newtheorem{lemma}{Lemma}

\newcounter{example}[section]
\newenvironment{example}[1][]{\refstepcounter{example}}{\medskip}
\makeatletter
\usetikzlibrary{calc}
\newcommand{\DrawLine}{%
  \begin{tikzpicture}
  \path[use as bounding box] (0,0) -- (\linewidth,0);
  \draw[color=black,solid,line width=1.2pt]
        ($(-0.3,0.1)$)--($(12.2,0.1)$);
  \end{tikzpicture}%
  }   
\makeatother

\title{Nonlinear Optimal Guidance for Fixed-Time Impact on a Stationary Target}

\author{
Kun Wang, Zheng Chen*, Han Wang, Jun Li\\
  \texttt{z-chen@zju.edu.cn}
   }

\affil{School of Aeronautics and Astronautics, Zhejiang University, Hangzhou 310027, Zhejiang, China}

\begin{document}

\maketitle

\begin{abstract}
This paper is concerned with devising the nonlinear optimal guidance for intercepting a stationary target with a fixed impact time. According to  Pontryagin's Maximum Principle (PMP), some optimality conditions for the solutions of the nonlinear optimal interception problem are established, and the structure of the corresponding optimal control is presented.  By employing the optimality conditions, we formulate a parameterized system so that its solution space is the same as that of the nonlinear optimal interception problem. As a consequence, a simple propagation of the parameterized system, without using any optimization method, is sufficient to generate enough sampled data for the mapping from  current state and time-to-go to the optimal guidance command. By virtue of the universal approximation theorem, a feedforward neural network, trained by the generated data, is able to represent the mapping from current state and time-to-go to the optimal guidance command. Therefore, the trained network eventually can  generate fixed-impact-time nonlinear optimal guidance within a constant time. Finally, the developed nonlinear optimal guidance is exemplified and studied through simulations, showing that the nonlinear optimal guidance law performs better than existing interception guidance laws. 
\end{abstract}

\section{ Introduction}
\label{intro}

Optimal guidance acts as an indispensable role during a homing missile approaches to a desired target because it determines the maneuvering commands for the missile  while optimizing a defined performance \cite{doi:10.2514/1.G006191}. By linearizing the engagement around a nominal collision course, the linearized/linear optimal guidance laws have been well devised since the 1960s \cite{Bryson:65,Willems:68,Cottrell:71,Anderson:79,Anderson:81}. The most popular example of optimal guidance, in the linearized setting, is probably the conventional Proportional Navigation (PN). Though PN was initially derived from physical intuition, it has been proven to be optimal with a coefficient of $3$ \cite{Kreindler:73,Bryson:75}.   
It is well known that the final impact time cannot be controlled if PN is used to guide an interceptor to a target. From a practical point of view, it is sometimes required to control or fix the final impact time. The scenario of salvo attack is a typical example requiring fixed-impact-time guidance because   multiple missiles should be coordinated to arrive a single target simultaneously \cite{doi:10.2514/1.40136}. Another example is that a missile intercepts a time-sensitive target, for which the missile should be guided to the target by a fixed duration of time. Due to the practical importance of controlling impact time, many scholars have devoted significant efforts to deriving fixed-impact-time guidance laws for interception problems. 

A popular way for devising fixed-time guidance laws is to adjust the navigation gain of PN according to the difference between desired time-to-go and estimated time-to-go. For this reason, some methods for estimating the time-to-go of PN have been developed in the literature; see, e.g., \cite{DAVID1991Time, 1996Closed,2015Accurate}. Based on the estimation of time-to-go, Jeon I.-S. {\it et al.} devised an Impact-Time-Control Guidance (ITCG) in the seminal paper \cite{Jeon:2006} in 2006.  The ITCG was further generalized in \cite{jeon2016impact} for intercepting moving targets based on the closed-form solution of PN.  Through deriving the analytical solution for the time-to-go in the nonlinear setting, a guidance law with time-varying guidance gain based on pure PN to control the impact time was proposed by Cho and Kim in \cite{2016Modified}. A generalized formulation for the ITCG in \cite{Jeon:2006,jeon2016impact}  was developed by He and Lee \cite{doi:10.2514/1.G003343} through studying optimal error dynamics; using the similar method, an extension to 3-dimensional interception scenarios was presented later on in  \cite{doi:10.2514/1.G003971}.

Advanced control theories have also been employed to devise fixed-impact-time guidance laws. Based on sliding mode methods, some guidance laws to control the impact time were developed \cite{2012Impact,2015Impact,Cho2016Nonsingular,2018Nonsingular,kumar2015impact}. Saleem and Ratnoo \cite{doi:10.2514/1.G001349} proposed an exact closed-form impact time expression based on Lyapunov theory, and the control command was implemented by varying a control parameter during the interception. The Lyapunov stability theory was employed by Kim {\it et al.} \cite{7126168}  to devise $2$- and $3$-dimensional nonlinear guidance laws with impact time constraints. In addition to the constraint on impact time, more complex constraints such as impact angle and field of view have also been taken into account in designing guidance laws; see, e.g., \cite{lee2007guidance,kim2013augmented}. Recently, taking  into account constraints on field of view, a guidance strategy for intercepting a stationary target with a fixed impact time was established by using a barrier Lyapunov function in \cite{9409726}.
 Through analyzing the geometric properties of the interception trajectory, a circular impact-time guidance was developed by Tsalik and Shima \cite{tsalik2019circular}; this method also applies to generating guidance laws for intercepting moving targets with both impact-angle and look-angle constraints.

According to the review in the preceding paragraphs, the optimal guidance laws with constraint on impact time  were  derived either  using  advanced control theories or based on variants of  PN. Though using the theories of sliding mode control and Lyapunov stability to devise guidance laws enables to satisfy complex constraints,  the optimality of control effort cannot be guaranteed. In the linearized setting, PN-based fixed-impact-time guidance laws can ensure the optimality to some extent. However, once the deviation from collision triangle is relatively large, the performance of PN-based guidance laws  can be improved if the nonlinear kinematics is considered \cite{chen2019nonlinear,Guelman:84,Lu:06,Jeon:10}.

The optimal guidance with nonlinear kinematics is generally dubbed Nonlinear Optimal Guidance (NOG), and it becomes an active research topic in the field of guidance in recent decades. In essence, the NOG is determined by the solutions of a class of nonlinear optimal control problems. As guidance commands should be generated onboard or in real time, it is usually required to devise a closed-form solution for the optimal control problems. The first attempt to devising closed-form solution for the NOG was probably done by  Guelman and Shinar \cite{Guelman:84}, indicating that the NOG is determined by a zero of  three nonlinear equations with three unknown variables. Numerical or optimization methods have also been developed in order for generating NOG in real time. For instance, Sim, Leng, and Subramaniam \cite{Sim:2000} proposed to combine a genetic algorithm with  a  shooting method to generate optimal solutions. Liu, Shen, and Lu  \cite{Liu:17}  developed a successive convex optimization approach  to generate NOG by transforming the nonlinear optimal control problems to a nonlinear programming problem. More recently, a parameterization method was developed by Chen and Shima \cite{chen2019nonlinear} so that the NOG can be generated efficiently and robustly by finding the zeros of a real-valued function. 

The final impact time is assumed free in all the papers cited in the previous paragraph. Recently, a quadratic method was proposed  by   Merkulov,   Weiss, and   Shima in \cite{merkulov2021minimum} to approximate the nonlinear kinematics in the vicinity of the initial Line-Of-Sight (LOS). This allowed establishing  a semi-analytic solution to the Fixed-Time NOG (FTNOG for abbreviation), and it presented an important attempt to devising FTNOG. Whereas, to the authors' best knowledge, the study on FTNOG without any approximation on nonlinear kinematics is rare to see in the literature. 


The current paper, as a continuing effort of \cite{chen2019nonlinear} where the final time is free, aims to develop a real-time method for generating the FTNOG through employing the optimality conditions and geometric properties of the solutions of nonlinear optimal interception problem. First of all, a new optimality condition is established by analyzing the necessary conditions from Pontryagin's Maximum Principle (PMP). This new optimality condition indicates that the look angle along an optimal trajectory cannot be equal to zero or $\pi$ (cf. Lemma \ref{LE:cos_sigma=1}). This further implies that the optimal control cannot change its sign more than once (cf. Lemma \ref{LE:control_twice}).

According to the necessary conditions from the PMP as well as the new established optimality condition, it is found in the paper that the optimal solutions are determined by two bounded scalars. Using the two bounded scalars, a parameterized system is formulated so that its  solution space is the same as the space of the solutions of the fixed-time nonlinear optimal interception problem. As a result, a simple propagation of the parameterized system gives rise to enough sampled data for the mapping from current state and time-to-go to the optimal control command. Then, in virtue of the universal approximation theorem \cite{hornik1989multilayer}, a neural network, trained by the sampled data, is sufficient to generate the optimal control command  within a constant time. 


As a matter of fact, because of the powerful capability of approximating highly nonlinear mappings, neural networks have been applied to many fields in aerospace engineering, such as orbital transfer \cite{cheng2018real, peng2018artificial, doi:10.2514/1.G005254}, powered descent guidance \cite{you2021learning,cheng2020real}, trajectory optimization \cite{shi2021onboard}, etc. Using a neural network for guidance relies heavily on  generating a large number of optimal trajectories as the training data set. In the literature, the data set is usually generated  via either indirect method, often combined with homotopy (continuation) method \cite{doi:10.2514/1.G002357,shi2021onboard} or direct method \cite{sanchez2016learning}. However, both methods suffer the issue of convergence \cite{chen2019nonlinear}. Another issue is that the solution generated by either indirect or direct method cannot be guaranteed to be optimal; see, e.g., \cite{chen20161,chen2021second}. If non-optimal solutions exist in the data set, the trained network cannot be used to generate optimal control command. 
In contrary, the developments in the paper, without using any optimization methods, allow to generate the data set by a simple propagation, as shall be shown by Procedure  \ref{algo1}.  In addition, thanks to the new optimality condition in Lemma \ref{LE:cos_sigma=1}, all the mappings from current state and time-to-go to the optimal control command generated by Procedure \ref{algo1} are related to the global optimal solutions. This ensures the trained network to generate global optimal control command, as illustrated by the numerical simulations in Section \ref{SE:Numerical}.

The paper is organized as follows: The nonlinear optimal control problem for the FTNOG is formulated  in Section \ref{SE:problem}. In Section \ref{SE:properties}, optimality conditions are derived according to the PMP, and a parameterized system for the optimal solutions is formulated; in addition, the geometric properties of the optimal control command is studied. 
In Section \ref{SE:Real}, the procedure for generating the optimal control command   by training a neural network   is presented. In Section \ref{SE:Numerical}, some numerical examples are presented, demonstrating and verifying the developments of the paper. This paper finally concludes by Section \ref{SE:conclusions}. 

\section{Problem Formulation}\label{SE:problem}

\begin{figure}[!htp]
\begin{center}
\includegraphics[trim=1cm 1.5cm 1cm 1cm, clip=true, width=3.5in]{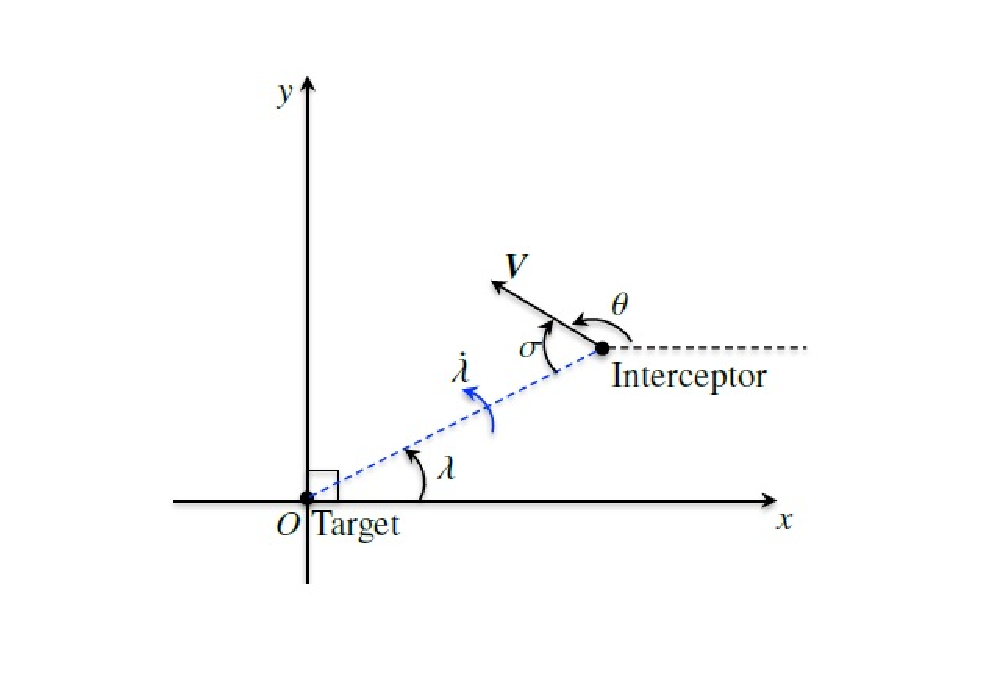}
\caption{The 2-dimensional intercept geometry.}\label{Fig:frame}
\end{center}
\end{figure}

Consider the 2-dimensional interception geometry on the horizontal plane in Fig.~\ref{Fig:frame} with a non-moving target. The origin of the  frame $Oxy$ is located at the target, the $x$ axis points to the East, and the $y$ axis points to the North.
Throughout the paper, we denote by $(x,y)\in\mathbb{R}^2$ the position of the interceptor in the frame $Oxy$, and by $\theta\in[0,2\pi]$  the angle between the $x$ axis and the velocity vector $\boldsymbol{V}$ of the interceptor. The angle $\theta$ is usually called the heading angle of the interceptor and it is positive when measured counterclockwise.
By normalizing the magnitude of the velocity $\boldsymbol{V}$ of the interceptor to one, the nonlinear kinematics of the interceptor  is
 represented by \cite{Lu:06}
\begin{align}
\begin{cases}
\dot{x}(t) =  \cos \theta(t),\\
\dot{y}(t) =  \sin \theta(t),\\
\dot{\theta}(t) = u(t),
\end{cases}
\label{EQ:Sigma}
\end{align}
where $t>0$ is the time, the dot  denotes the differentiation with respect to time, and $u\in\mathbb{R}$ is the control parameter which represents the lateral acceleration. Then, generating the FTNOG requires addressing the following Optimal Control Problem (OCP).
\begin{problem}\label{problem:OIP}
Given an initial condition $(x_0,y_0,\theta_0)\in \mathbb{R}^2\times [0,2\pi]$ and a fixed final time $t_f>0$,  the OCP consists of steering the system in Eq.~(\ref{EQ:Sigma}) by a measurable control $u(\cdot)$ on $[0,t_f]$ from the initial state $(x_0,y_0,\theta_0)$  to the final point $(0,0)$ with the final heading angle left free such that the control effort
$$ J = \int_0^{t_{f}}\frac{1}{2} u^2(t) \mathrm{d}t$$
is minimized. 
\end{problem}

In order to find the solution of the OCP within a constant time, the ideal way is about devising an analytical formula for the optimal solution. However, because of the nonlinear property in Eq.~(\ref{EQ:Sigma}), it is hardly possible to obtain an analytical formula; see \cite{Guelman:84,chen2019nonlinear}. In the following two sections, we shall show how to develop a real-time method for the optimal solutions by combining the properties of the optimal trajectories with artificial neural network. 


\section{Characterizations of the Optimal Solutions}\label{SE:properties}

In this section, some optimality conditions will be first established according to the PMP. Then, by using these optimality conditions, the optimal solutions are parameterized by two bounded scalars. In addition, the geometric properties of the optimal control will be presented. For simplicity of presentation, important results and claims are written in lemmas and the proofs are postponed to appendices.

\subsection{Optimality Conditions}\label{SE:Pontryagin}

Denote by $p_x$, $p_y$, and $p_{\theta}$ the co-state variables of $x$, $y$, and $\theta$, respectively. Then, the Hamiltonian for the OCP is expressed as
\[H= p_x \cos \theta + p_y \sin \theta + p_{\theta} u + p^0  \frac{1}{2} u^2,\]
where $p^0$ is a non positive scalar. In view of \cite[Remark 2]{chen2019nonlinear}, we have that $p^0$ is negative. For any negative $p^0$, the quadruple $(p_x,p_y,p_{\theta},p^0)$ can be normalized so that $p^0 = -1$. Thus, we shall
consider $p^0 = -1$ in the remainder of the paper.

According to the PMP in \cite{Pontryagin}, we have    
\begin{align}
\begin{cases}
\dot{p}_x(t) = -\frac{\partial H}{\partial x} =  \ 0,\\
\dot{p}_y(t) = -\frac{\partial H}{\partial y} =  \ 0,\\
\dot{p}_{\theta}(t) = -\frac{\partial H}{\partial \theta} = \ p_x(t) \sin \theta(t) - p_y(t) \cos \theta(t),
\end{cases}
\label{EQ:dot_p_theta}
\end{align}
and
\begin{align}
\frac{\partial H}{\partial u} = 0.\label{EQ:dH/du}
\end{align}
Explicitly rewriting Eq.~(\ref{EQ:dH/du}) leads to
\begin{align}
 u(t) =  p_\theta(t),\ \ t\in[0,t_f].
\label{EQ:u=p_theta}
\end{align}
Since the final heading angle is free, the transversality condition implies
\begin{align}
p_{\theta}(t_f) = 0.
\label{EQ:Transversality1}
\end{align}
 In view of Eq.~(\ref{EQ:dot_p_theta}), we have that $p_x$ and $p_y$ are constants along an optimal trajectory. Taking into account the final boundary condition $(x(t_f),y(t_f)) = (0,0)$ and   Eq.~(\ref{EQ:Transversality1}), we can integrate the third equation of Eq.~(\ref{EQ:dot_p_theta}) to yield
\begin{align}
p_{\theta}(t) = p_x y(t) - p_y x(t),\ t\in[0,t_f].
\label{EQ:p_theta}
\end{align}

For brevity, the triple $(x(t),y(t),\theta(t))$ for $t\in [0,t_f]$ is said to be an {\it extremal trajectory} if all the necessary conditions in Eqs.~(\ref{EQ:dot_p_theta}--\ref{EQ:p_theta}) are met. According to \cite{chen20161,chen2021second}, any extremal trajectory cannot be guaranteed to be  optimal unless sufficient optimality conditions are satisfied. By the following lemmas, we establish a new optimality condition. 
\begin{lemma}\label{LE:cos_sigma=1}
Given any extremal trajectory $(x(t),y(t),\theta(t))$ for $t\in[0,t_f]$, if there exists a time $\bar{t}\in (0,t_f)$ so that  the velocity vector $[\cos \theta(\bar{t}),\sin \theta(\bar{t})]$ is collinear with the LOS, i.e., 
\begin{align}
\frac{|
x(\bar{t})\cos \theta(\bar{t}) + y(\bar{t})\sin \theta(\bar{t})|}{
\sqrt{x(\bar{t})^2 + y(\bar{t})^2}} = 1
\label{EQ:cos_sigma=1}
\end{align}
then the extremal trajectory $(x(\cdot),y(\cdot),\theta(\cdot))$ on $[0,t_f]$ is not optimal.
\end{lemma}
The proof of this lemma is postponed to Appendix \ref{Appendix:A}.  Thanks to Lemma \ref{LE:cos_sigma=1}, we immediately have the following result on the structure of optimal control.
\begin{lemma}\label{LE:control_twice}
Given any extremal trajectory $(x(t),y(t),\theta(t))$ for $t\in[0,t_f]$, let $u(t)$ for $t\in [0,t_f]$ be the corresponding extremal control. If there are two different times $t_1\in (0,t_f)$ and $t_2\in (0,t_f)$ so that $u(t_1) = u(t_2) = 0$, then the extremal trajectory is not optimal. 
\end{lemma}
The proof of this lemma is postponed to Appendix \ref{Appendix:A}. 

Notice that Lemma \ref{LE:cos_sigma=1} presents a new optimality condition, and Lemma \ref{LE:control_twice} indicates that the optimal control cannot change its sign more than once. 

\subsection{Parameterized System}

Define the constant $\alpha\geq 0$ as the norm of the vector $[p_x,p_y]$, i.e., 
$$\alpha \coloneqq \sqrt{p_x^2 + p_y^2},$$
and define the constant $\beta\in[-\pi,\pi]$ as the solution of 
\begin{align}
p_x = \alpha \cos \beta\ \text{and}\ p_y =  \alpha \sin \beta.
\label{EQ:definition_alpha}
\end{align}
Then, we can rewrite Eq.~(\ref{EQ:p_theta}) as
\begin{align}
p_{\theta}(t) = \alpha[y(t) \cos \beta - x(t) \sin \beta].
\label{EQ:p_theta1}
\end{align}
Let us define a parametrized system
\begin{align}
\begin{cases}
\dot{X}(t) =  -\cos \Theta(t),\\
\dot{Y}(t) =  -\sin \Theta(t),\\
\dot{\Theta}(t) = -\alpha [Y(t)\cos \beta - X(t) \sin \beta],
\end{cases}
\label{EQ:Sigma1}
\end{align}
where $(X,Y)\in \mathbb{R}^2$ and $\Theta \in [0,2\pi]$. 
For the sake of notational simplicity, let the triple 
$$(X(t,\alpha,\beta),Y(t,\alpha,\beta),\Theta(t,\alpha,\beta))\in \mathbb{R}^2\times [0,2\pi]$$
 for $t\in [0,t_{f}]$ be the solution of the  $(\alpha,\beta)$-parameterized system  in Eq.~(\ref{EQ:Sigma1}) with the initial condition of $(0,0,0)$. In other words, we have  
$$(X(0,\alpha,\beta),Y(0,\alpha,\beta),\Theta(0,\alpha,\beta)) = (0,0,0)$$ for any $\alpha>0$ and $\beta\in [-\pi,\pi]$. In the next section, we shall show that the solution space of the parameterized system is the same as the solution space of OCP in the polar coordinate.

\subsection{Optimal Control Command}\label{SE:Policy}

Let $\lambda \in [0,2\pi]$ be the angle between the LOS and the $x$ axis, as shown in Fig.~\ref{Fig:frame}. This angle is positive if measured counterclockwise. Then, we have that the nonlinear kinematics in the polar frame is represented by
\begin{align}
\begin{cases}
\dot{r}(t) = \cos [ \theta(t) - \lambda(t)],\\
\dot{\lambda}(t) = \sin [ \theta(t) - \lambda(t)]/r(t),\\
\dot{\theta}(t) = u
\end{cases}
\label{EQ:Sigma_r}
\end{align}
where $r\geq 0$ is the Euclidean distance between the interceptor and the target, i.e., 
\begin{align}
r = \sqrt{x^2 + y^2}
\end{align}
Let $\sigma\in [-\pi,\pi]$ be the angle between the LOS and velocity vector $\boldsymbol{V}$,   as shown in Fig.~\ref{Fig:frame}. This angle is called the look angle, which is positive when measured clockwise. As $\sigma = \pi +  \lambda - \theta$,  we can further write the system in Eq.~(\ref{EQ:Sigma_r}) as 
\begin{align}
\begin{cases}
\dot{r}(t) = - \cos \sigma(t),\\
\dot{\sigma}(t) = \frac{ \sin  \sigma(t)}{r(t)} - u(t),\\
\end{cases}
\label{EQ:Sigma_s}
\end{align}
Since the kinematics in Eq.~(\ref{EQ:Sigma_s}) is equivalent to that in Eq.~(\ref{EQ:Sigma}) before $r$ approaches to zero, it follows that the optimal control command can be either determined by the state $(x,y,\theta)$ in Cartesian frame or by the state $(r,\sigma)$ in the polar frame. 

Let us denote by $C(r_c,\sigma_c,t_g)$ the optimal control command at the current state $(r_c,\sigma_c)$ with a feasible time-to-go of $t_g>0$.  This means that, given any optimal trajectory $(r(t),\sigma(t))$ for $t\in [0,t_f]$, if $u(t)$ is the corresponding optimal control, then we have $C(r(t), \sigma(t),t_f - t) = u(t)$ for $t\in [0,t_f]$.  With the definition of optimal control command, it is clear that the  real-time generation of optimal solutions is equivalent to computing the value of the mapping $(r_c,\sigma_c,t_g)\mapsto C(r_c,\sigma_c,t_g)$  in real time. By the following lemma, we present a symmetric property for the optimal control command. 
\begin{lemma}\label{LE:opposite}
Given an current state $(r_c,\sigma_c)$ in polar frame and a feasible time-to-go $t_{g}>0$, we have 
$$C(r_c,\sigma_c,t_g) = - C(r_c,-\sigma_c,t_g).$$ 
\end{lemma}
The proof of this lemma is postponed to Appendix \ref{Appendix:A}. 


Note that
\begin{align}
\cos \sigma(t) = 
\frac{
x( {t})\cos \theta( {t}) + y( {t})\sin \theta( {t})}{
\sqrt{x( {t})^2 + y( {t})^2}}
\label{EQ:sigma_t}
\end{align}
According Eq.~(\ref{EQ:sigma_t}) and Eq.~(\ref{EQ:cos_sigma=1}) in Lemma \ref{LE:cos_sigma=1}, given any optimal trajectory $(r(t),\sigma(t))$ for $t\in [0,t_f]$ in the polar frame, we have that the look angle $\sigma(\cdot)$ on the open interval $(0,t_f)$ cannot be equal to either $0$ or $\pi$. This further indicates that along any optimal trajectory we have either $\sigma(t)\in (0,\pi)$ or $\sigma(t)\in (-\pi,0)$ holds for all $t\in (0,t_f)$. In view of Lemma \ref{LE:opposite},  given an current state $(r_c,\sigma_c)$ and a feasible time-to-to $t_g>0$, if $\sigma_c\in [-\pi,0]$, we can use $-C(r_c,-\sigma_c,t_g)$ as the optimal control command for the current state $(r_c,-\sigma_c)$ with a time-to-go of $t_g$.  As a result, in order to generate the optimal control command in real time, it amounts to developing a real-time method for computing the value of $C(r_c,\sigma_c,t_g)$ with $\sigma_c$ taking values only in the semi-open interval $(0,\pi]$. 

For notational simplicity, we denote by
\begin{align}
\mathcal{F}_P = \{(r_c,\sigma_c,t_g) | r_c>0,\sigma_c \in [0,\pi], t_g > 0 \}
\end{align}
the set of current state $(r_c,\sigma_c)$ associated with feasible time-to-go $t_g$.  In other words, given any $(r_c,\sigma_c,t_g)\in \mathcal{F}_P$, there exists an optimal trajectory $(r(t),\sigma(t))$ for $t\in [0,t_g]$  so that $(r(0),\sigma(0)) = (r_c,\sigma_c)$ and $r(t_g) = 0$. In the next section, we shall show how the optimality conditions and properties established in this section endow the capability of generating the value of $C(r_c,\sigma_c,t_g)$ for any $(r_c,\sigma_c,t_g)\in \mathcal{F}_P$ in real time via artificial neural networks. 


\section{Real-Time Solution for the OCP}\label{SE:Real}



According to the universal approximation theorem \cite{hornik1989multilayer}, if we use a large number of sampled data of the relationship between $(r_c,\sigma_c,t_g)$ and $C(r_c,\sigma_c,t_g)$  to train a feedforward neural network, then the trained network is able to  accurately represent the mapping $(r_c,\sigma_c,t_g)\mapsto C(r_c,\sigma_c,t_g)$  for any  $(r_c,\sigma_c,t_g ) \in \mathcal{F}_P$. Since the output of a feedforward neural network is a composition of some linear mappings of input vector, it follows that  for any input vector $(r_c,\sigma_c,t_g ) \in \mathcal{F}_P$, the trained neural network can  generate the corresponding optimal control command  $ C(r_c,\sigma_c,t_g)$ within a constant time, which paves a way for generating solutions of OCPs in real time. This is why many scholars have used neural networks to develop real-time method for OCPs in aerospace engineering, as stated in Section \ref{intro}.

A prerequisite of training a  neural network to generate the optimal control command  is to generate the data set for the relationship between  $(r_c,\sigma_c,t_g)$ and $C(r_c,\sigma_c,t_g)$. A straightforward way to generate the data set is to use optimization methods (including indirect methods and direct methods) to solve the corresponding OCP, as was done in, e.g., \cite{doi:10.2514/1.G002357,sanchez2016learning}. However, both indirect and direct methods suffer the issues of convergence. While homotopy method can be combined with indirect method to improve the convergence  \cite{shi2020deep}, it is well known that homotopy method may fail to converge as well \cite{caillau2012differential}. On the other hand, even if an optimization method converges, the generated solution cannot be guaranteed to be at least locally optimal  \cite{chen2021second,chen20161}. If non-optimal solutions are included in the data set, the trained network cannot be ensured to generate optimal control command. In the following paragraphs, we shall show how to employ the developments in Section \ref{SE:properties} to generate the data set for the optimal solution without any optimization. 



Set $\Lambda = (0,+\infty)\times [0,\pi]$. For any $t>0$ and $(\alpha,\beta)\in \Lambda$, let us set
\begin{align}
\begin{cases}
R(t,\alpha,\beta)=   \sqrt{X^2(t,\alpha,\beta)+Y^2(t,\alpha,\beta)}\\
\Sigma (t,\alpha,\beta)=    \arccos \Big\{ - \frac{X(t,\alpha,\beta) \cos[ \Theta(t,\alpha,\beta)] + Y(t,\alpha,\beta) \sin[ \Theta(t,\alpha,\beta)] }{\sqrt{X^2(t,\alpha,\beta) + Y^2(t,\alpha,\beta)  }} \Big \}\\
U(t,\alpha,\beta) =  \alpha [Y(t,\alpha,\beta) \cos \beta -X(t,\alpha,\beta) \sin \beta ] 
\end{cases}
\label{EQ:U}
\end{align}
and
\begin{align}
T(\alpha,\beta) =  \min \{t> 0 | \cos [\Sigma(t,\alpha,\beta)] = 1  \} 
\end{align} 
By the definition of $T(\alpha,\beta)$, we have that the parameterized trajectory $(X(t,\alpha,\beta),Y(t,\alpha,\beta),\Theta(t,\alpha,\beta))$ for $t\in [0,T(\alpha,\beta)]$ not only satisfies the necessary conditions from PMP but also satisfies the optimality condition in Lemma \ref{LE:cos_sigma=1}. 
Then, we have the following result.
\begin{lemma}\label{LE:h=u}
Given any $(r_c,\sigma_c,t_g)\in \mathcal{F}_P$, there exists $(\alpha,\beta)\in \Lambda$ so that
\begin{align}
\begin{cases}
t_g \in   (0,T(\alpha,\beta)]\\
r_c =  R(t_g,\alpha,\beta)\\
\sigma_c =   \Sigma(t_g,\alpha,\beta)\\
C(r_c,\sigma_c,t_g) =  U(t_g,\alpha,\beta)
\end{cases}
\label{EQ:LE:h=u}
\end{align}
\end{lemma} 
The proof of this lemma is postponed to Appendix \ref{Appendix:B}. 

\begin{lemma}\label{LE:r_sigma_in_F}
Given any $(\alpha,\beta)\in \Lambda$, if $t\in [0,T(\alpha,\beta)]$, then we have
\begin{align}
\begin{cases}
(R(t,\alpha,\beta),\Sigma(t,\alpha,\beta),t)\in \mathcal{F}_P\\
C(R(t,\alpha,\beta),\Sigma(t,\alpha,\beta),t)= U(t,\alpha,\beta)
\end{cases}
\label{EQ:LE:r_sigma_in_F}
\end{align}
\end{lemma}
The proof of this lemma is postponed to Appendix \ref{Appendix:B}.






Lemmas \ref{LE:h=u} and \ref{LE:r_sigma_in_F} imply that the space of the  parameterized solutions in Eq.~(\ref{EQ:U}) is the same as the space of optimal solutions in polar coordinate. 
Thus, we can use the mapping $(R(t,\alpha,\beta),\Sigma(t,\alpha,\beta),t)\mapsto U(t,\alpha,\beta)$ for $(\alpha,\beta)\in \Lambda$ and $t\in (0,T(\alpha,\beta)]$ to generate the data set for the relationship between $(r_c,\sigma_c,t_g)$ and $ C(r_c,\sigma_c,t_g)$.
For each $(\alpha,\beta)\in \Lambda$ and $t\in [0,T(\alpha,\beta)]$,  the value of the mapping $(R(t,\alpha,\beta),\Sigma(t,\alpha,\beta),t)\mapsto U(t,\alpha,\beta)$ can be readily obtained by propagating the parameterized system in Eq.~(\ref{EQ:Sigma1}) with the initial condition $(0,0,0)$ and by using the transformation in Eq.~(\ref{EQ:U}). Note that the value of $T(\alpha,\beta)$ can be too large, resulting in the difficulty of  propagating Eq.~(\ref{EQ:Sigma1}) to generate the data set. To address this difficulty, let us set
\begin{align}
\hat{T}(\alpha,\beta) = \min\big\{\bar{T}, T(\alpha,\beta) \big\}
\label{EQ:hat_T}
\end{align}
where $\bar{T}>0$ is a finite number. Then,  the procedure for generating the sampled values for the relationship between $(r_c,\sigma_c,t_g)$ and $ C(r_c,\sigma_c,t_g)$ is summarized in Procedure~\ref{algo1} (from step 1 to step 7). Note that all the sampled data   are finally included in the set $\mathcal{D}$. 
\begin{table}
\usetikzlibrary{calc}
\begin{center}
	\begin{tcolorbox}[
		colframe=blue!25,
		colback=gray!20,
		coltitle=blue!20!black,  
		fonttitle=\bfseries,
		width = 0.8\columnwidth,
		boxrule=0pt,
		top=1pt,
		bottom=-0.5pt,
		enhanced,
		segmentation style={solid, black,line width=1.2pt},
		segmentation code={\draw[black,solid,line width=1.2pt]($(segmentation.west)+(0.2,0)$)--($(segmentation.east)+(-0.2,0)$);}]
		\begin{example}
		\label{algo1}
		\end{example}
\textbf{Procedure 1: Real-Time Generation of Optimal Solutions via  Neural Networks}
		\tcblower
\begin{itemize}
\item[1.] Let $\bar{\alpha}$, $h$, and $\bar{T}$ be positive numbers, and let $N_i$ and $N_j$ be positive integers. Set $i=1$, $j=1$, and $\mathcal{D}=\varnothing$.
\item[2.] If $i\leq N_i$, set $\alpha = i \bar{\alpha}/N_i $ and go to step 3; otherwise, go to step 8.
\item[3.] If $j\leq N_j$, set $\beta = j \pi/N_j$ and go to step 4; otherwise, set $i=i+1$ and go to step 2.
\item[4.] Propagate the parameterized system in Eq.~(\ref{EQ:Sigma1}) with the initial condition $(0,0,0)$ and use the transformation in Eq.~(\ref{EQ:U}) to generate the trajectory $(R(t,\alpha,\beta),\Sigma(t,\alpha,\beta))$ and the corresponding control $U(t,\alpha,\beta)$ for $t\in [0,\hat{T}(\alpha,\beta)]$. Set $t=0$ and go to step 5.
\item[5.] If $t+ h \leq \hat{T}(\alpha,\beta)$, set $t = t + h$ and go to step 6; otherwise, set $j=j+1$  and go to step 3.
\item[6.] Set 
\begin{align}
r=&\ R(t,\alpha,\beta)\nonumber\\
\sigma =&\ \Sigma(t,\alpha,\beta)\nonumber\\
u =&\ U(t,\alpha,\beta)\nonumber
\end{align}
and go to step 7.
\item[7.] Set $\mathcal{D} = \mathcal{D}\cup \{[r,\sigma,t,u]\}$, and go to step 5.
\item[8.] Use the data set $\mathcal{D}$ to train a neural network to represent the optimal control command $(r_c,\sigma_c,t_g)\mapsto C(r_c,\sigma_c,t_g)$. 
\end{itemize}
\DrawLine
	\end{tcolorbox}
\end{center}
\end{table}
In view of the universal approximation theorem   established in  \cite{hornik1989multilayer}, a simple feedforward   neural network, trained by the data in $\mathcal{D}$, is enough to represent the mapping $(r_c,\sigma_c,t_g)\mapsto C(r_c,\sigma_c,t_g)$, as shown by step 8 in Procedure~\ref{algo1}.

\begin{figure}[!htp]
\centering
\subcaptionbox{Current state measured in polar frame\label{guidance2}}{\includegraphics[height = .19\linewidth]{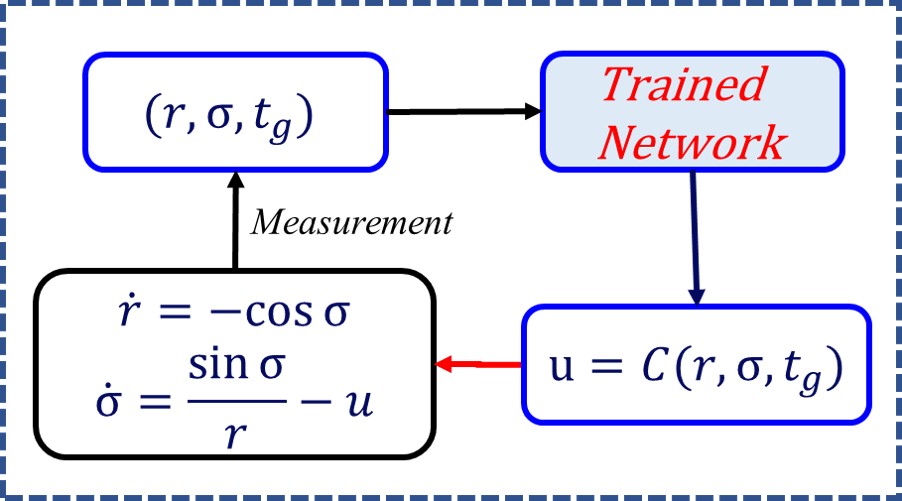}}
~~~~~~
\subcaptionbox{Current state measured in frame $Oxy$\label{guidance1}}{\includegraphics[height = .19\linewidth]{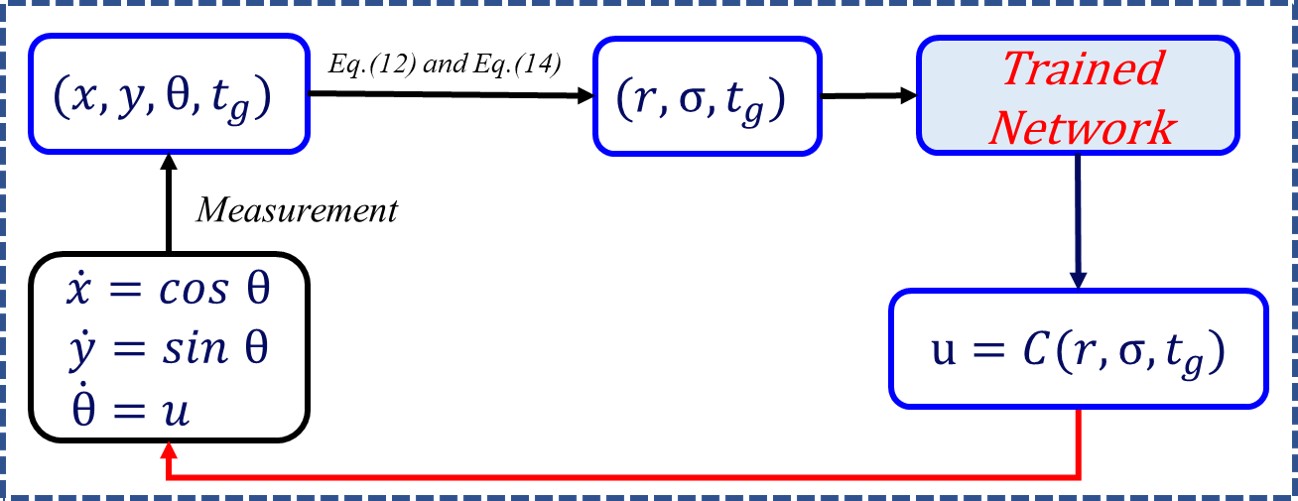}}
\caption{Diagrams for closed-loop NOG.}
\label{Fig:diagram_guidance}
\end{figure}

Because $\bar{T}$ is a finite number, it follows that the data set $\mathcal{D}$ does not contain all the data for the relationship between $(r_c,\sigma_c,t_g)$ and $ C(r_c,\sigma_c,t_g)$. By the following lemma, we shall show that even if $\bar{T}$ is a small number, the trained neural network in Procedure \ref{algo1}, is enough to represent the  mapping $(r_c,\sigma_c,t_g)\mapsto C(r_c,\sigma_c,t_g)$ for any $(r_c,\sigma_c,t_g)\in \mathcal{F}_P$. 
\begin{lemma}\label{LE:Norm}
Let $(\bar{r}(t),\bar{\sigma}(t))$ for $t\in [0,\bar{t}_g]$ be an optimal trajectory of an interceptor whose speed is $\bar{V}>0$, and let $\bar{u}({t})$ for $t\in [0,\bar{t}_g]$  be the corresponding optimal control. Then, for any $t_g \in (0,\bar{t}_g)$, we have 
\begin{align}
\bar{u}({t}) = \frac{\bar{V} t_g}{\bar{t}_g} C \left(\frac{\bar{r}({t}) t_g}{\bar{V}\bar{t}_g},\bar{\sigma}({t}),{t_g} - {t}\frac{t_g}{\bar{t}_g}\right)
\label{EQ:LE:Norm}
\end{align}
\end{lemma}
The proof of this lemma is postponed to Appendix \ref{Appendix:B}. Notice that, given any current state $(\bar{r}_c,\bar{\sigma}_c)$ with $\bar{\sigma}_c \in [0,\pi]$ and a feasible time-to-go $\bar{t}_g>0$ for an  interceptor with any positive speed $\bar{V}>0$, we can choose a small time-to-go $t_g > 0$ and  a pair $(\alpha,\beta)\in \Lambda$ so that 
\begin{align}
\frac{\bar{r}_c t_g}{\bar{V}\bar{t}_g} =&\  R(t_g,\alpha,\beta)\nonumber\\
\bar{\sigma}_c =&\  \Sigma(t_g,\alpha,\beta)\nonumber\\
t_g \in &\  (0,\hat{T}(\alpha,\beta)]\nonumber
\end{align}
Therefore, even if $\bar{T}$ in Eq.~(\ref{EQ:hat_T}) is a small number, we can combine the trained network in Procedure \ref{algo1} and the transformation in Eq.~(\ref{EQ:LE:Norm}) to generate the optimal control command for an interceptor with any positive speed and any positive time-to-go.  The trained network can be employed in the closed-loop guidance, as shown in Fig.~\ref{guidance2}. 
It should be noted that if the measured state is $(x,y,\theta)$ in Cartesian coordinate instead of $(r,\sigma)$ in polar coordinate, the trained neural network can also be used by a simple transformation, as shown in Fig.~\ref{guidance1}. In the next section,  the performance of the developments will be numerically presented and verified by applying to practical examples.




\section{Numerical Simulations}\label{SE:Numerical}


Set the parameters $\bar{\alpha}$, $N_i$, $N_j$, $\bar{T}$, and $h$ in Procedure~\ref{algo1} as $10$, $100$, $100$, $10$, and $0.005$, respectively. With such parameters, the number of elements in the data set $\mathcal{D}$ generated in Procedure~\ref{algo1} is no more than $4.59\times 10^6$.  By using the data set $\mathcal{D}$, a neural network with three hidden layers, each of which contains 30 neurons, is trained to represent the optimal control command $C(r_c,\sigma_c,t_g)$ for $(r_c,\sigma_c,t_g)\in \mathcal{F}_P$. Hyperbolic tangent function is employed for all the hidden layers, and the training is terminated when the mean squared error between the predicted values and the values in $\mathcal{D}$ is less than $10^{-4}$.                                                                                                                                                                                       

Given any $(r_c,\sigma_c,t_g)\in \mathcal{F}_P$ as input, the trained network takes around $0.0935$ ms to generate an output using the Python library on a laptop with Intel Core i5-9300H CPU @2.40 GHz. This indicates the trained network can be used to generate the optimal control command with a frequency of more than $10$ KHz on the same computation platform. Notice that such a computational period is more than enough for the guidance systems of usual missiles, and the computational time can be reduced further if the trained model is implemented in a customized environment for onboard applications \cite{shi2020deep}. In the following two subsections, we shall show the performance of network-based FTNOG, by comparing with PN, ITCG \cite{Jeon:2006}, and optimal solutions from optimization methods.


\subsection{Performance of the FTNOG}\label{SE:Performance}


\paragraph{Case A: Comparisons with Different Impact Times}

Consider the interception problem in \cite{merkulov2021minimum}, where the initial position of the interceptor, denoted by $I$ hereafter, is at the origin $(0, 0)$ m with an initial heading angle of $60$ deg, and the speed is $500$ m/s. The target, denoted by $T$ hereafter, is located at $(10000, 0)$ m. Note that the target position can be transformed from $(0, 0)$ in Problem \ref{problem:OIP} to any desired position through a simple coordinate  transformation. We consider to combine the normalization in Lemma \ref{LE:Norm} and the trained network to generate FTNOG for the interceptor  with different desired impact times: $25$ sec, $30$ sec, $40$ sec, and $50$ sec. 

\begin{figure}[!htp]
\centering
\includegraphics[width = 1\linewidth]{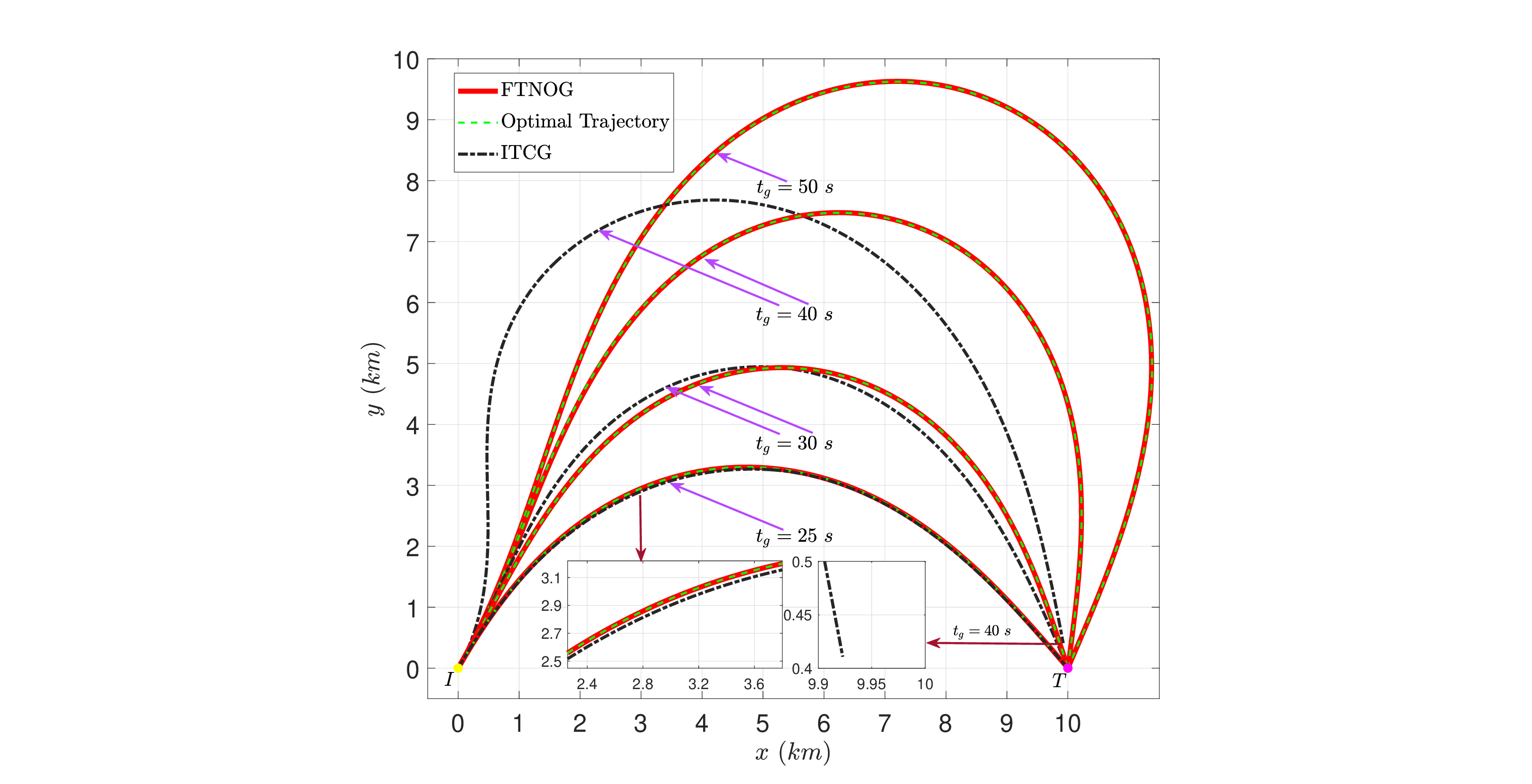}
\caption{Case A: Comparisons of trajectories related to different guidance laws with different impact times.}
\label{Fig:trajectories_35}
\end{figure}

The trajectories related to FTNOG are presented by the red solid curves in Fig.~\ref{Fig:trajectories_35}. Optimal solutions are presented by the green dashed curves in Fig.~\ref{Fig:trajectories_35}. It is clearly seen that the FTNOG-related trajectories coincide with the optimal solutions. The trajectories generated by the ITCG developed in \cite{Jeon:2006} are presented by black dotted-dashed curves in Fig.~\ref{Fig:trajectories_35}. 
We can see that the trajectories related to the ITCG diverge a lot from the optimal solutions, especially   when the desired impact times are large. It should be noted here that the ITCG cannot be used when the impact time is up to 50 sec because the error of time-to-go estimation becomes negative \cite{Jeon:2006}. Therefore, the ITCG-related trajectory with $t_g = 50$ sec is not presented in Fig.~\ref{Fig:trajectories_35}. Regarding the case that the initially desired impact time is  $40$ sec, the miss distance from the interceptor to the target for the ITCG is $417.99$ m when the time-to-go decreases to zero; see the scaled plot in Fig.~\ref{Fig:trajectories_35}. In contrary, the miss distance from the interceptor to the target for the FTNOG is quite close to zero for any desired impact times. 

\begin{figure}[!htp]
\centering
\begin{subfigure}[t]{7cm}
\centering
\includegraphics[width = 8cm]{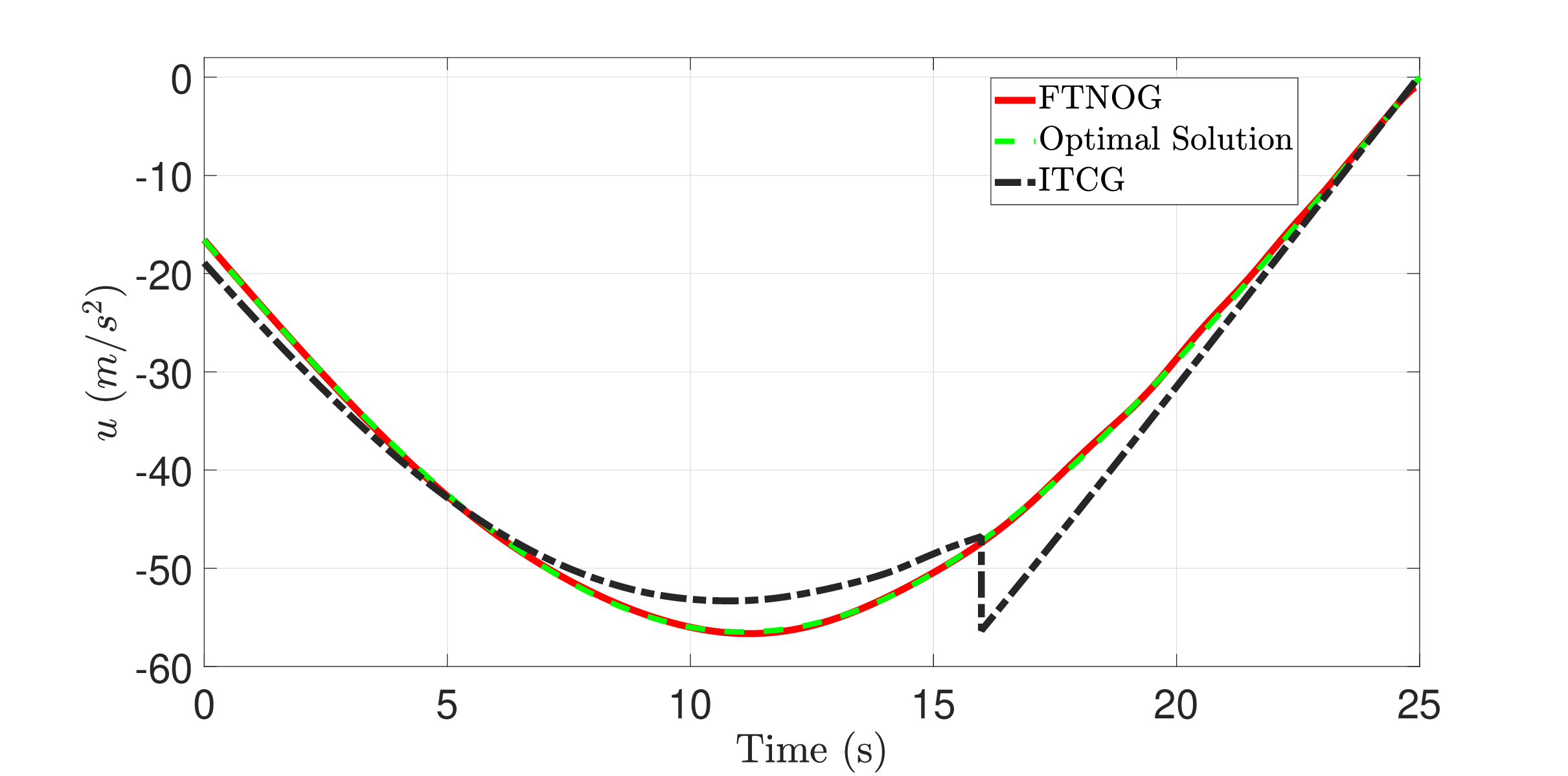}
\caption{ Desired impact time of 25 s}
\label{Fig:25_control}
\end{subfigure}
~~~~~
\begin{subfigure}[t]{7cm}
\centering
\includegraphics[width = 8cm]{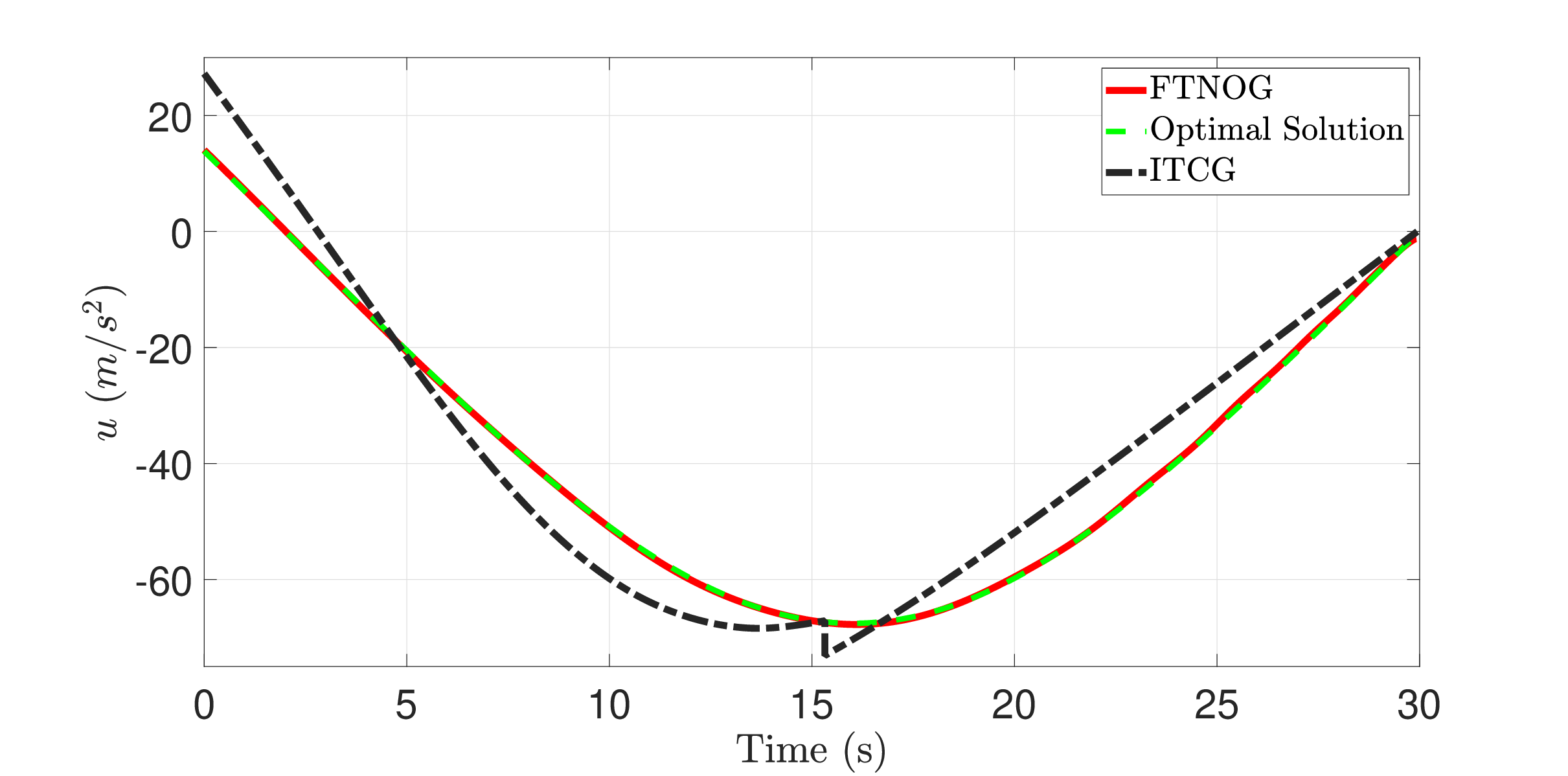}
\caption{Desired impact time of 30 s}
\label{Fig:30_control}
\end{subfigure}\\
\begin{subfigure}[t]{7cm}
\centering
\includegraphics[width = 8cm]{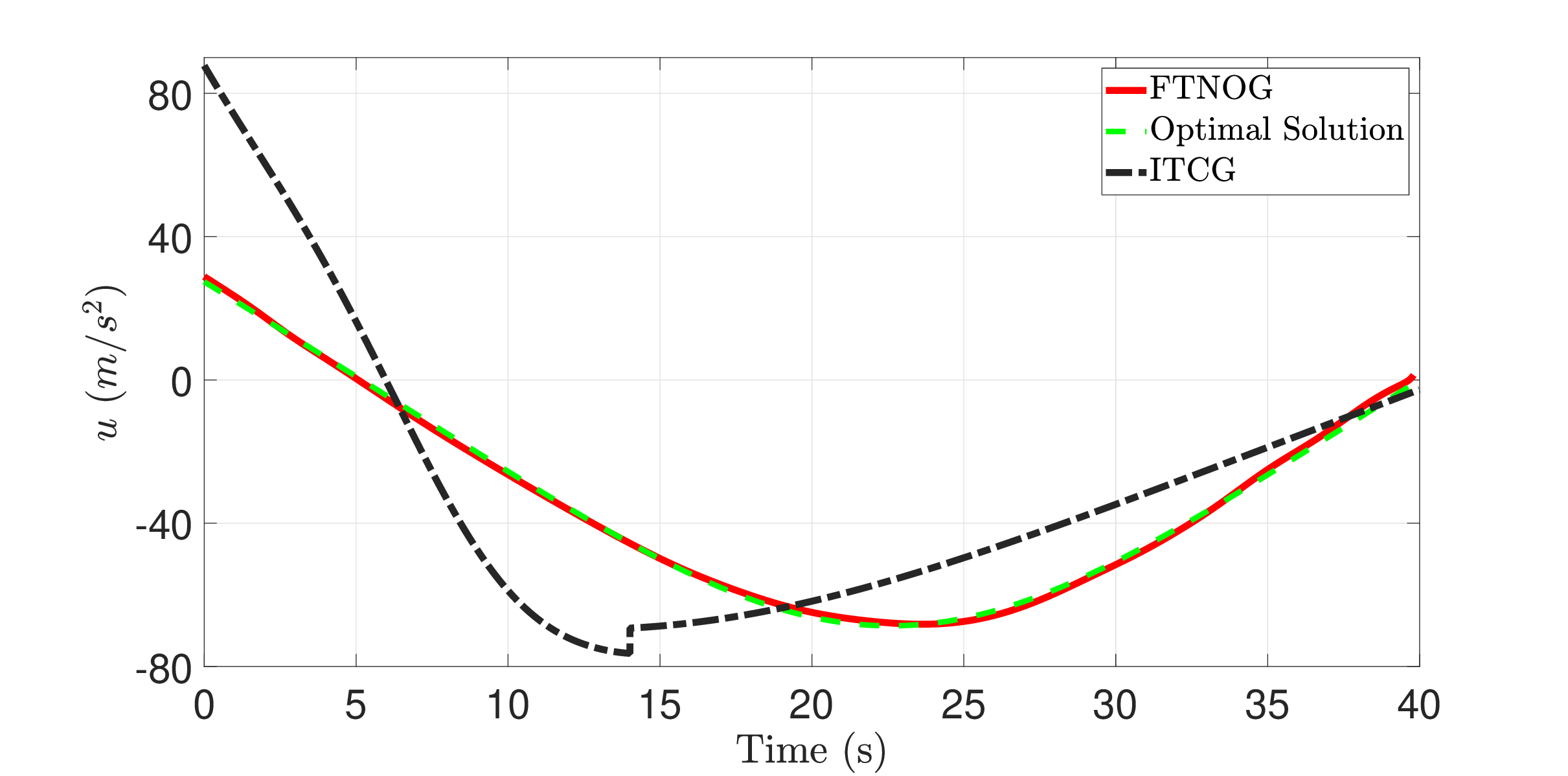}
\caption{Desired impact time of 40 s}
\label{Fig:40_control}
\end{subfigure}
~~~~~
\begin{subfigure}[t]{7cm}
\centering
\includegraphics[width = 8cm]{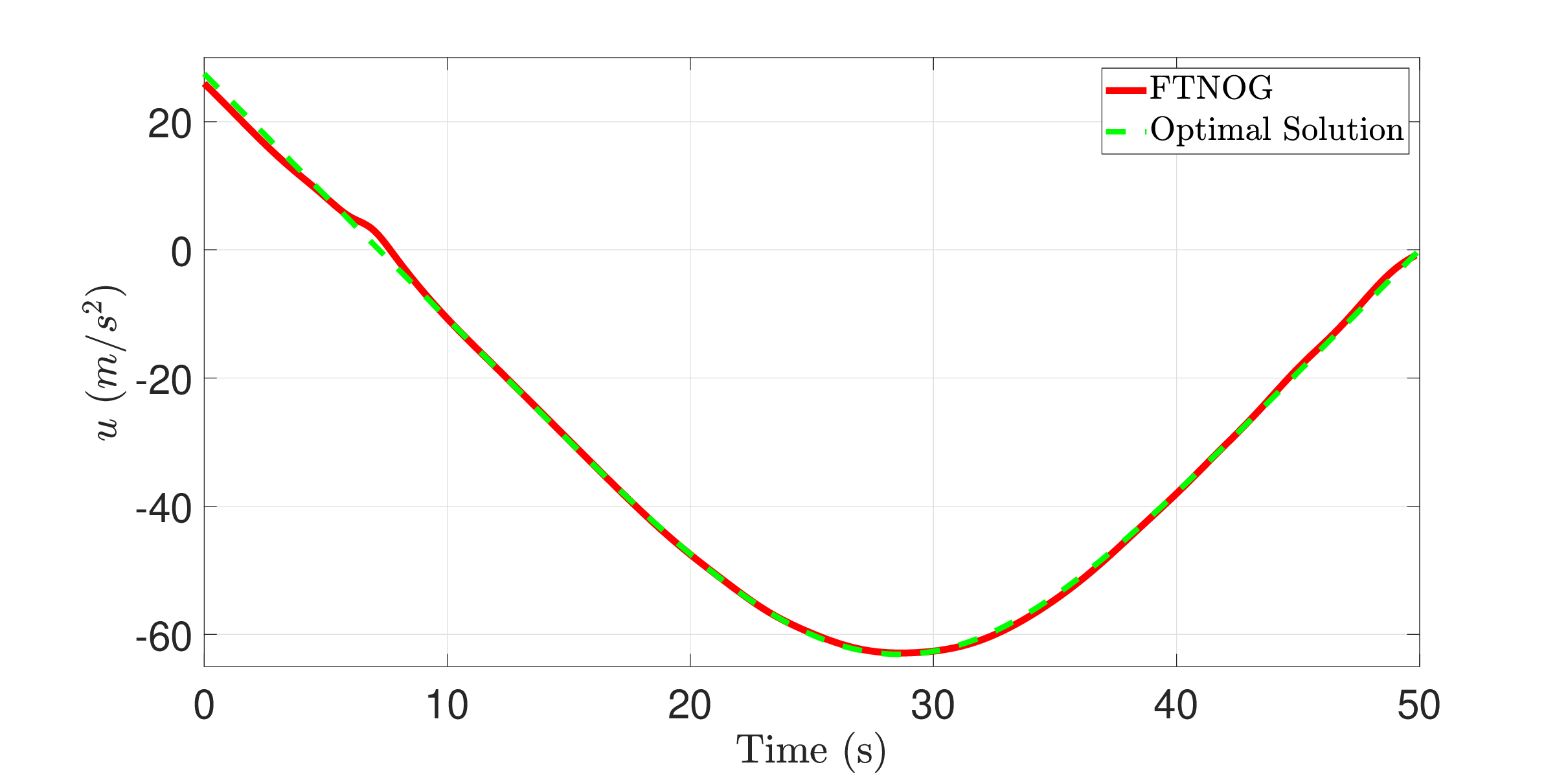}
\caption{Desired impact time of 50 s}
\label{Fig:50_control}
\end{subfigure}
\caption{Case A: Control profiles   with different impact times.}
\label{Fig:control_profile}
\end{figure}
\begin{figure}[!htp]
\centering
\includegraphics[width = 1\linewidth]{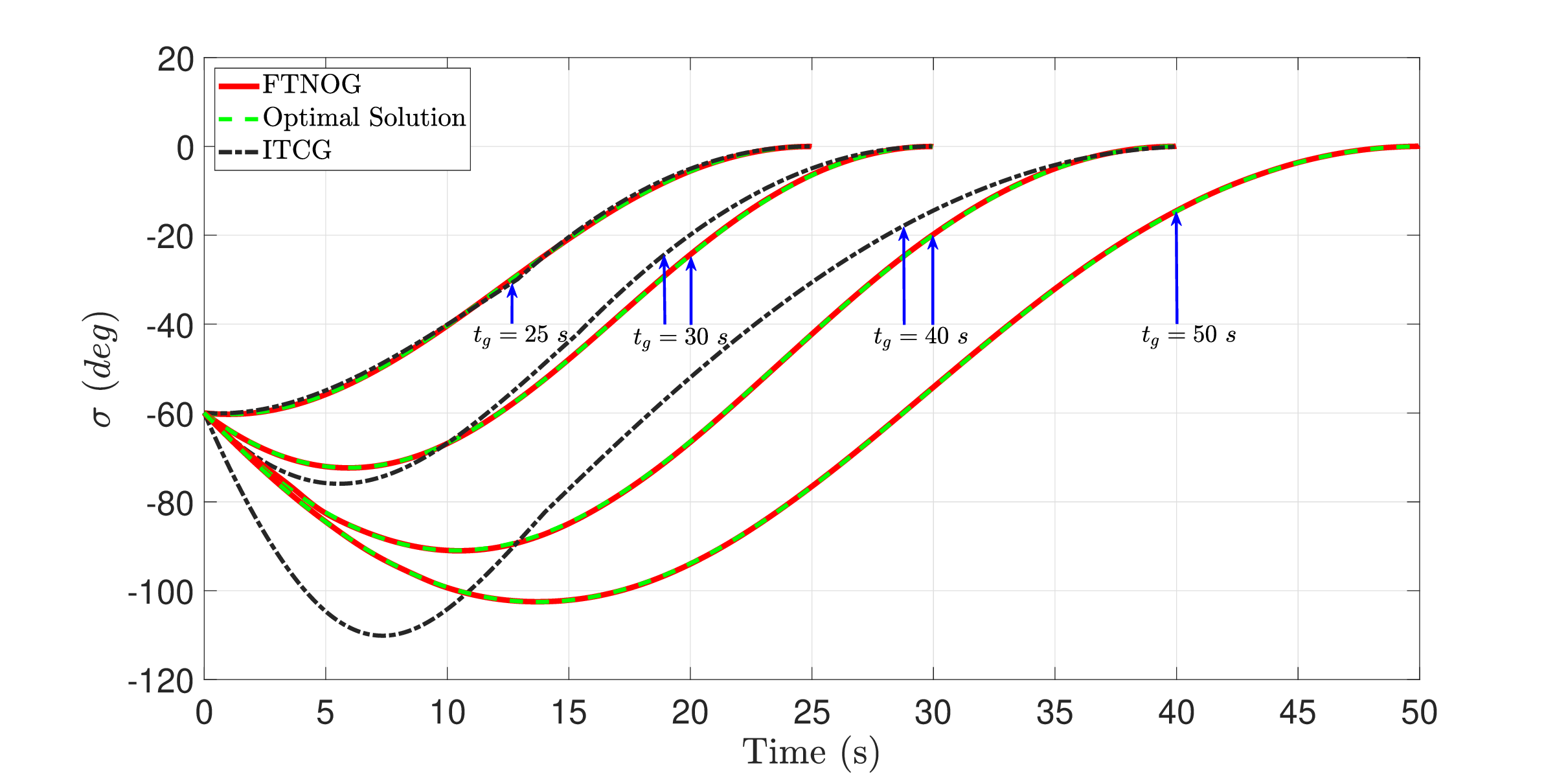}
\caption{Case A: Time history of look angle with different impact times.}
\label{Fig:Look angle_profile}
\end{figure}

 The profiles of corresponding controls and look angles are presented in Fig.~\ref{Fig:control_profile} and Fig.~\ref{Fig:Look angle_profile}, respectively.
Note that  the ITCG should be switched to PN with a navigation gain of 3 once the error of time-to-go estimation reaches zero \cite{Jeon:2006}. Thus, the controls are not continuous for the ITCG, as shown by the black dotted-dashed curves in Fig.~\ref{Fig:control_profile}. We can also see that the controls generated by FTNOG is a.e. the same as the optimal controls. It can be seen from 
Fig.~\ref{Fig:Look angle_profile} that the absolute values of look angles along optimal trajectories do not reach the value of $0$ deg or $180$ deg, as predicted by Lemma \ref{LE:cos_sigma=1}. 
Let $J_{ITCG}$, $J_{FTNOG}$, and $J^*$ be the control effort required by ITCG, FTNOG, and optimal control, respectively, and they are presented in Table.~\ref{Table:control_effort_25_50}, where the term ``NA'' indicates that the corresponding guidance law does not apply. It can be observed that the deviations of the control efforts between FTNOG and optimal solutions are very small. However, the difference of control efforts required by ITCG and by optimal control becomes increasingly large when the desired time-to-go increases. 
\begin{table}[!htp]
\centering
\caption{Case A: The values of control effort consumed by different guidance laws with different impact times.}
\begin{tabular}{cccc}
\hline
      & $J_{ITCG}$~($m^2/s^3$)       & $J_{FTNOG}$~($m^2/s^3$)       & $J^*$($m^2/s^3$)   \\ 
\hline
$t_g = 25~s$  & $2.1434\times10^4$    & $ 2.1371\times10^4$   & $ 2.1350\times10^4$   \\ 
\hline
$t_g = 30~s$ & $3.1057\times10^4$      &$3.0624\times10^4$  & $3.0563\times10^4$         \\ 
\hline
$t_g = 40~s$ & $1.3487\times10^5$     & $ 3.8912\times10^4$      & $3.8738\times10^4$                       \\
\hline
$t_g = 50~s$ &NA      & $3.9678\times10^4$ & $3.9625\times10^4$                          \\
\hline
\label{Table:control_effort_25_50}
\end{tabular}
\end{table}

\paragraph{Case B: Comparisons with Different Initial Heading Angles}

\begin{figure}[!htp]
\centering
\includegraphics[width = 1\linewidth]{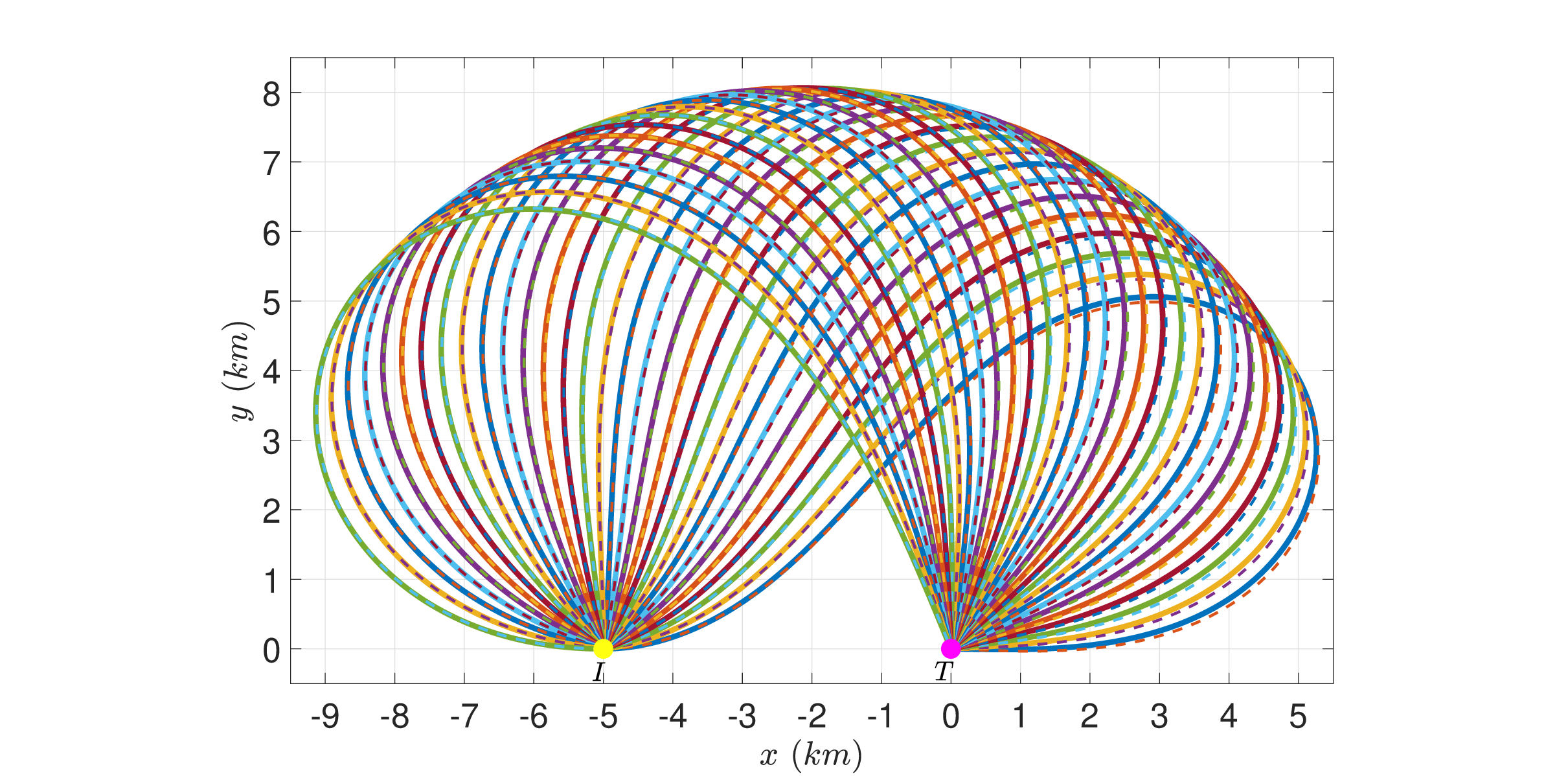}
\caption{Case B: Trajectories related to the FTNOG with different initial heading angles.}
\label{Fig:wide_applicability}
\end{figure}

For case B, we consider an interception scenario with different initial heading angles. The speed is set as $500$ m/s, and the desired impact time is set as $40$ sec. The initial position is set as $(-5000,0)$ m. The position of the target is   $(0,0)$ m. Then, the normalization in Lemma \ref{LE:Norm} and the trained network are combined to generate the FTNOG with different initial heading angles. The optimal trajectories (dashed curves) and FTNOG-related trajectories (solid curves) are presented in Fig.~\ref{Fig:wide_applicability}.

\begin{figure}[!htp]
\centering
\includegraphics[width = 1\linewidth]{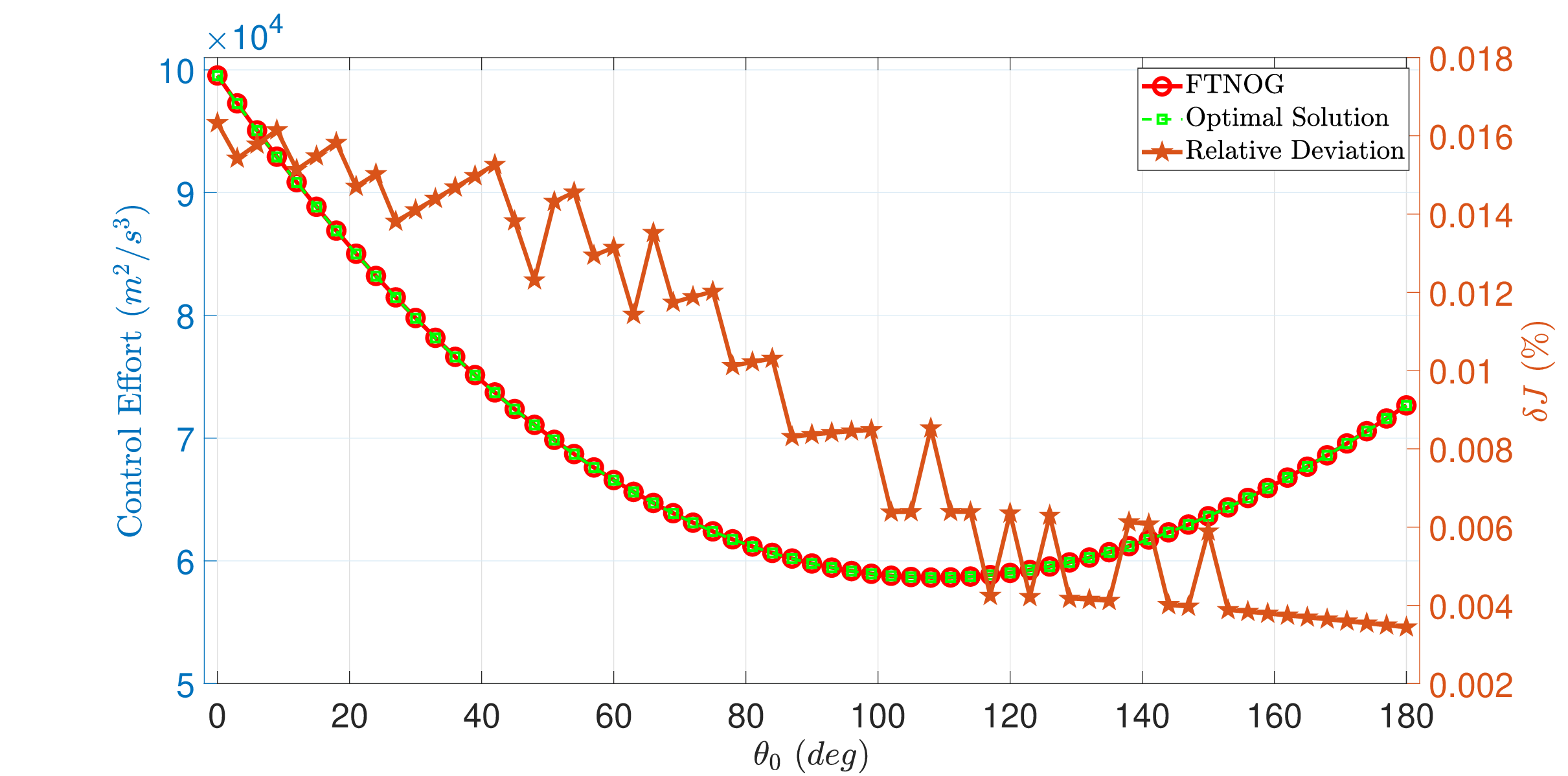}
\caption{Case B: Comparison of control efforts between the FTNOG and optimal solutions.}
\label{Fig:wide_applicability_cost}
\end{figure}

We can see from Fig.~\ref{Fig:wide_applicability} that if the initial heading angle is relatively large, the optimal trajectories almost coincide with the FTNOG-related trajectories. Note that if the initial heading angle approaches to zero, small deviations appear between optimal trajectories and FTNOG-related trajectories. However, the differences of control efforts are negligible even if the initial heading angle is close to zero, as shown by the control efforts by FTNOG and the control efforts by optimal control in Fig.~\ref{Fig:wide_applicability_cost}. The left vertical axis of Fig.~\ref{Fig:wide_applicability_cost} represents the values of control efforts, and the right vertical axis shows the relative deviation of $J_{FTNOG}$ from $J^*$, i.e.,
\begin{align}
\delta J = \frac{J_{FTNOG}-J^*}{J^*}\times100\% \nonumber
\end{align} 
According to Lemma \ref{LE:opposite}, the optimal trajectory with the initial condition of $(x_0,y_0,\theta_0)$ is symmetric to the optimal trajectory with  the initial condition of $(x_0,y_0,-\theta_0)$. Thus, only the trajectories with initial heading angles in $[0,180]$ deg are illustrated here. 

\begin{figure}[!htp]
\centering
\includegraphics[width = 1\linewidth]{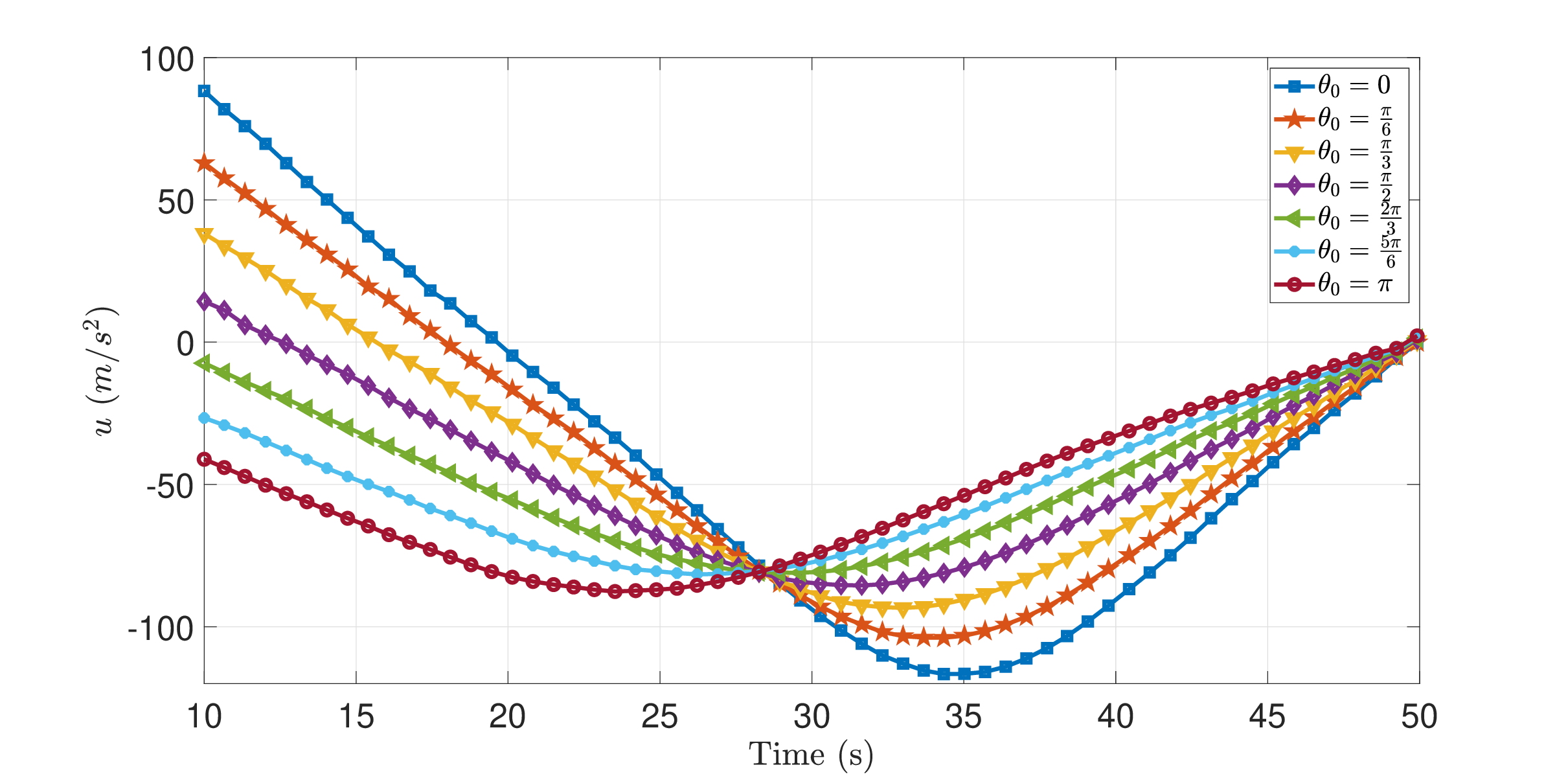}
\caption{Case B: Control profiles along trajectories related to the FTNOG with different initial heading angles.}
\label{Fig:wide_applicability_control}
\end{figure}

\begin{figure}[!htp]
\centering
\includegraphics[width = 1\linewidth]{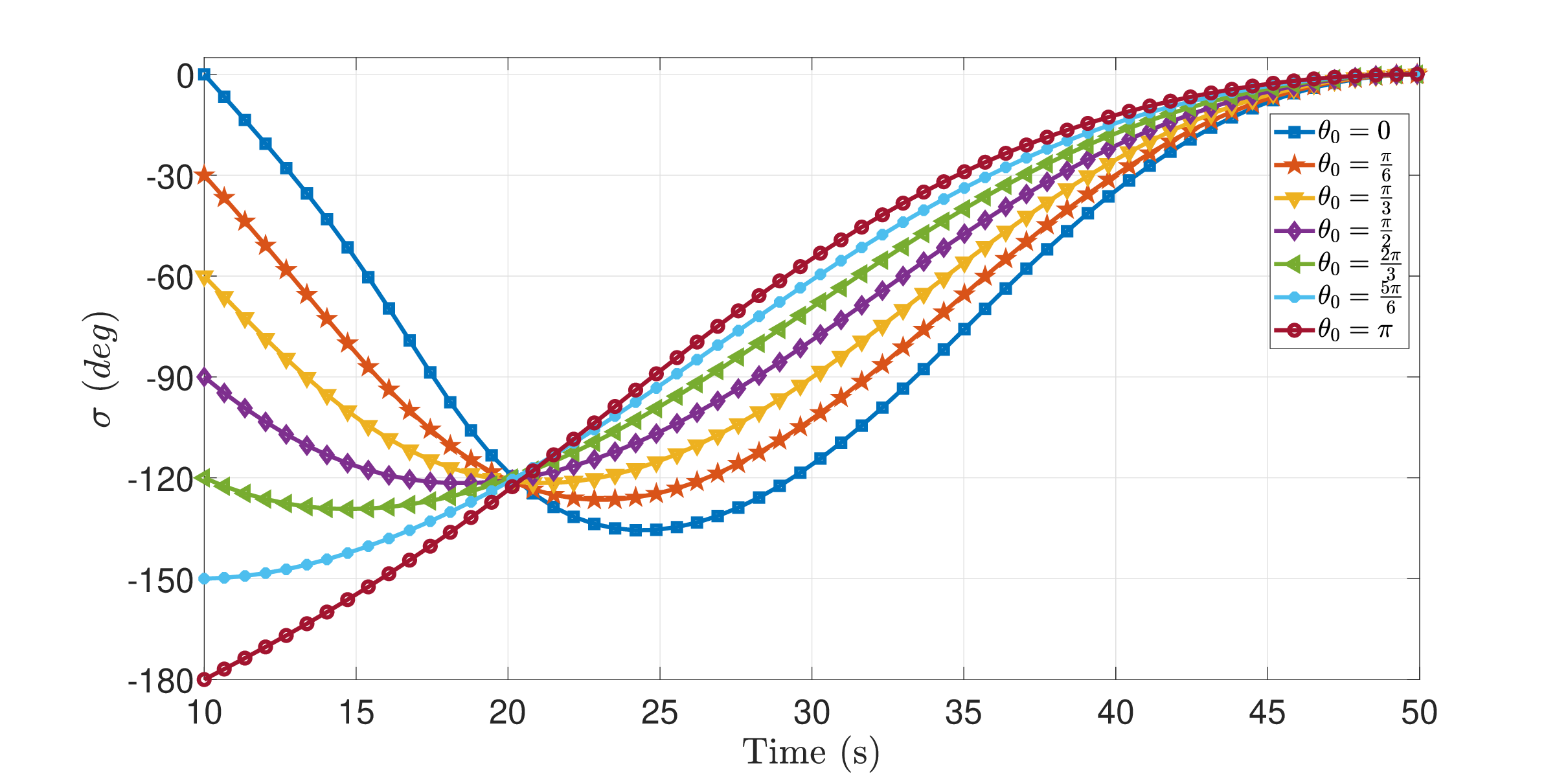}
\caption{Case B: Profiles of look angle along trajectories related to the FTNOG with different initial heading angles.}
\label{Fig:wide_applicability_sigma}
\end{figure}

Profiles of controls and look angles  along some FTNOG-related trajectories are shown in Fig.~\ref{Fig:wide_applicability_control} and Fig.~\ref{Fig:wide_applicability_sigma}, respectively. 
We can see from Fig.~\ref{Fig:wide_applicability_sigma} that the absolution values of look angles are in $[0,180]$ deg, and they never reach $0$ deg or $180$ deg with the exception of the endpoints, as predicted by the theoretical result in Lemma \ref{LE:cos_sigma=1}.

\begin{figure}[htbp]
 \centering 
\includegraphics[width = 1\linewidth]{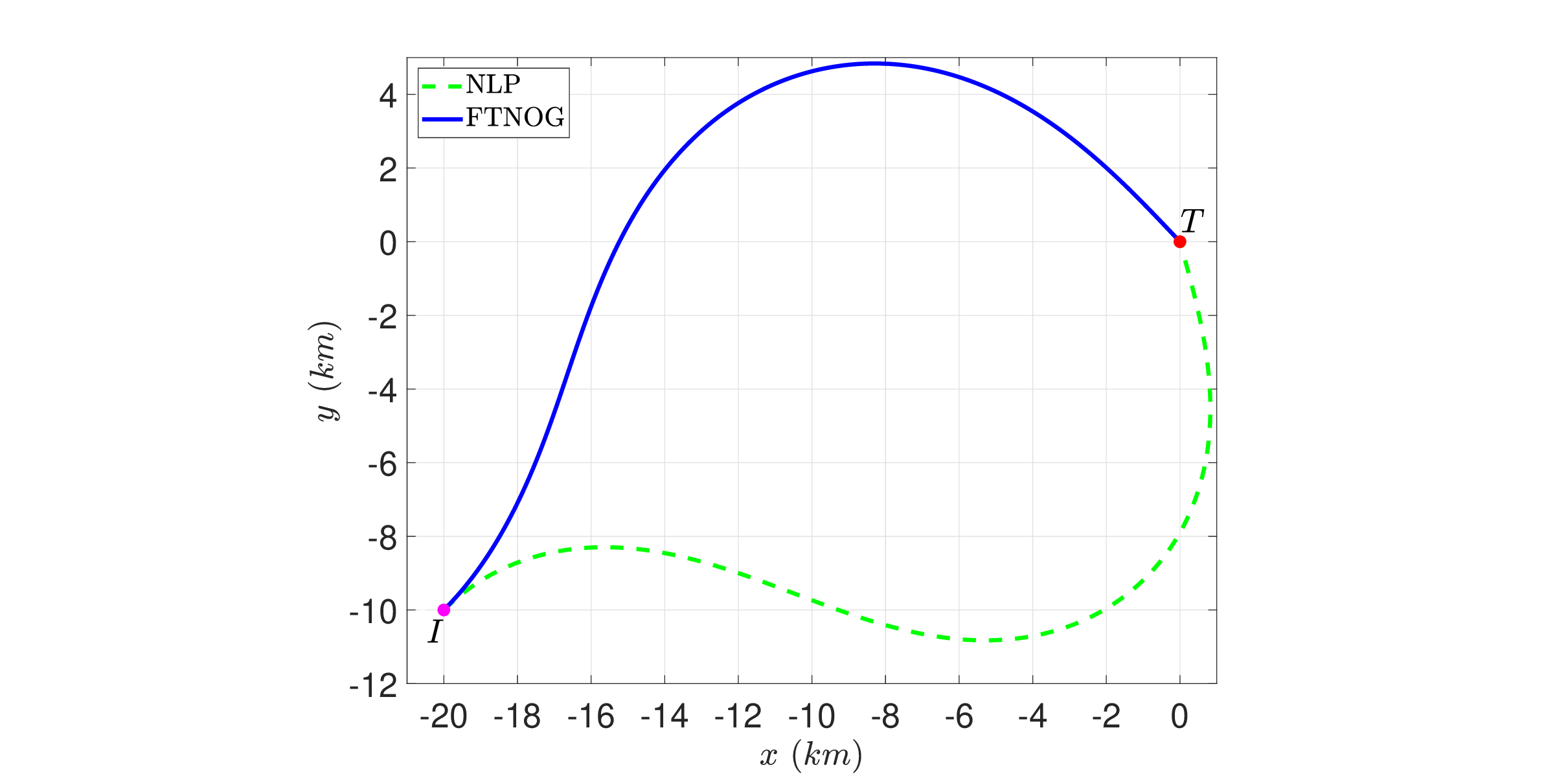}
\caption{Case C: Comparisons between trajectory related to FTNOG and trajectory from NLP.}
\label{Fig:GPOPS_DNN}
\end{figure}

\begin{figure}[htbp]
 \centering 
\includegraphics[width = 1\linewidth]{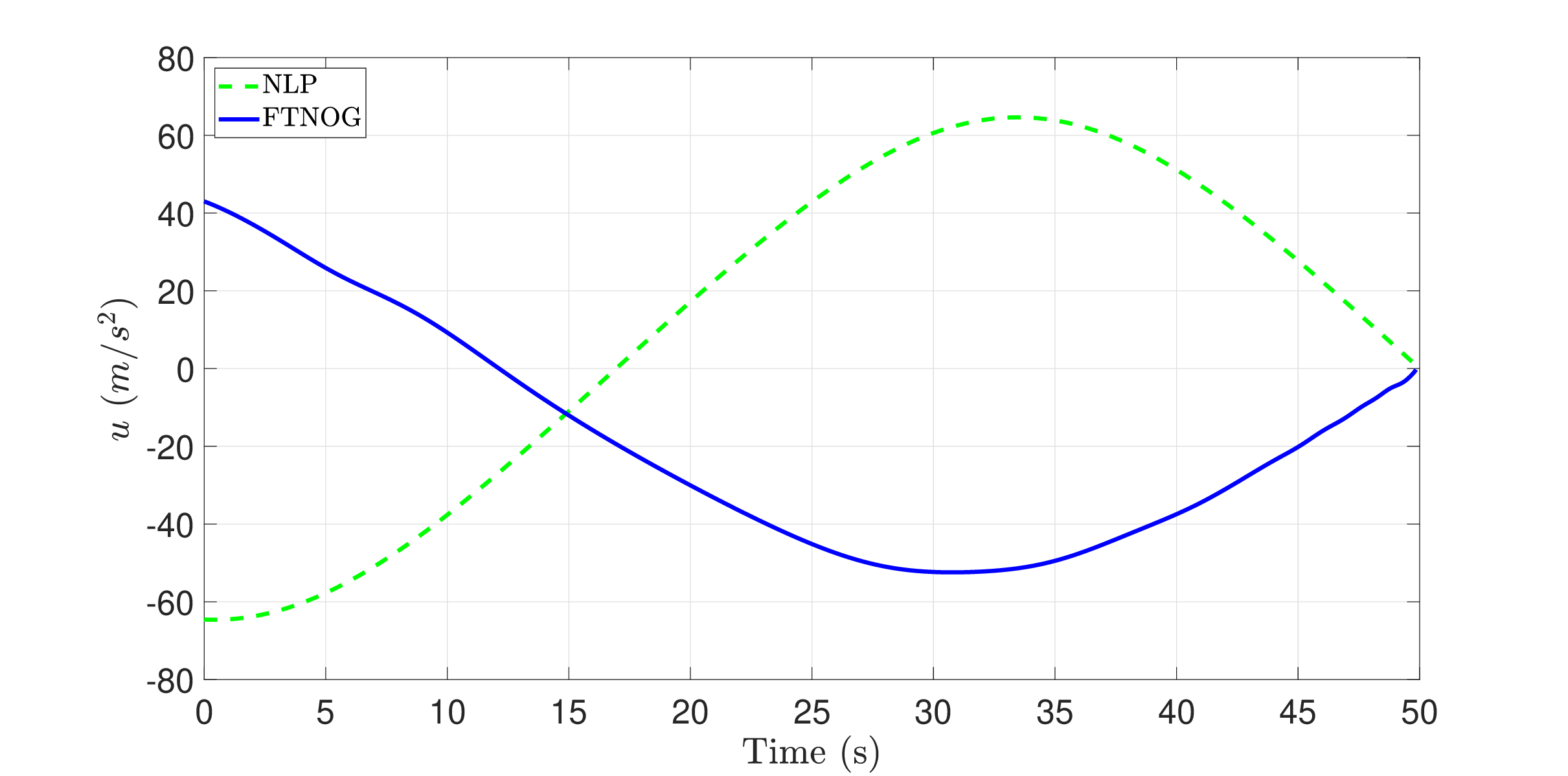}
\caption{Case C: Control profiles along the trajectory related to the FTNOG and the trajectory from NLP.}
\label{Fig:GPOPS_DNN1}
\end{figure}

\begin{figure}[htbp]
 \centering 
\includegraphics[width = 1\linewidth]{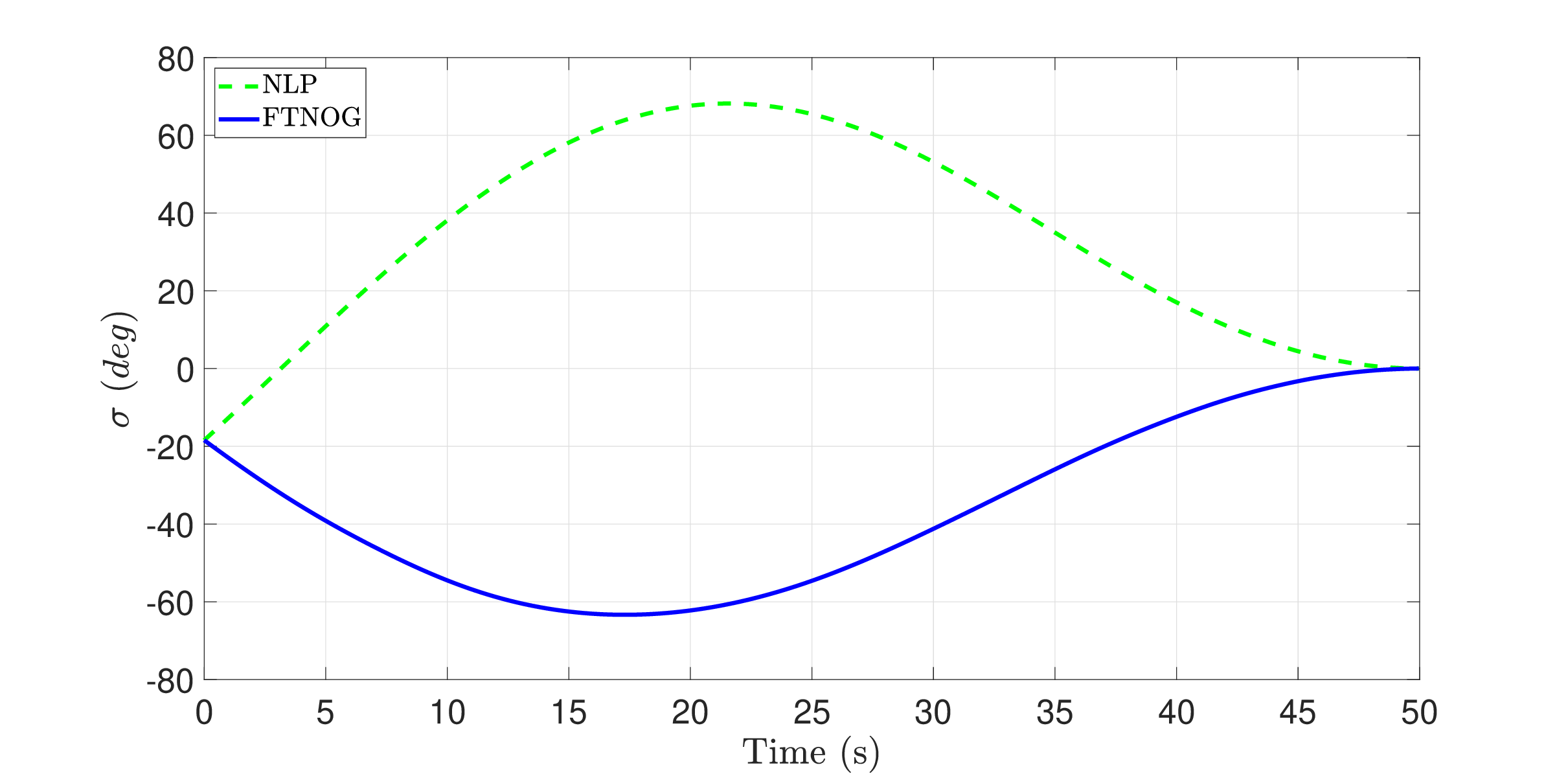}
\caption{Case C: Profiles of look angles along the trajectory related to the FTNOG and the trajectory from NLP.}
\label{Fig:GPOPS_DNN2}
\end{figure}

\paragraph{Case C: Capability of Generating Global Solutions}
In this paragraph, we consider to demonstrate the capability of the FTNOG to generate the global solutions. The initial state for the interceptor is set as
\begin{align}
 (x_0,y_0,\theta_0) = (-20000 ~\text{m},\ -10000 ~\text{m} ,\ \frac{\pi}{4} )\nonumber 
\end{align}
The speed is $600$ m/s, and the fixed impact time is set as $50$ sec. Then, the trained network, by combining with the normalization in Lemma \ref{LE:Norm}, is employed to generate the FTNOG, which finally yields a solution as shown by the solid curve in Fig.~\ref{Fig:GPOPS_DNN}. However, employing a Nonlinear Programming (NLP) method (the package of GPOPS-II with Radau orthogonal collocation method \cite{Gpops2} is used here) generates a totally different trajectory, as shown by the dashed curve in Fig.~\ref{Fig:GPOPS_DNN}.
The control profiles generated by FTNOG and NLP  are shown in Fig.~\ref{Fig:GPOPS_DNN1}. 
The control effort along the trajectory generated by NLP is $5.0572\times10^4$ $m^2/s^3$, while that related to FTNOG is only  $2.9158\times10^4$ $m^2/s^3$. Thus, the trajectory from NLP is not optimal. This happens because NLP methods usually find solutions by satisfying only  necessary conditions (e.g., Karush–Kuhn–Tucker conditions).  The profiles of look angles along the two trajectories are demonstrated in Fig.~{\ref{Fig:GPOPS_DNN2}}, 
from which we can see that there exists a point ( at around $3.15$ sec) along the trajectory generated by NLP so that the look angle is zero. This means that the solution from NLP does not meet the optimality condition in Lemma \ref{LE:cos_sigma=1}. Thanks to the theoretical developments in Section \ref{SE:properties}, all the non-optimal solutions are not included in the data set $\mathcal{D}$  in Procedure \ref{algo1}. Thus, the FTNOG generated by the trained network can lead to a global solution.

\subsection{Application to Salvo Attack Scenario}

\begin{figure}[!htp]
\centering
\includegraphics[width = 1\linewidth]{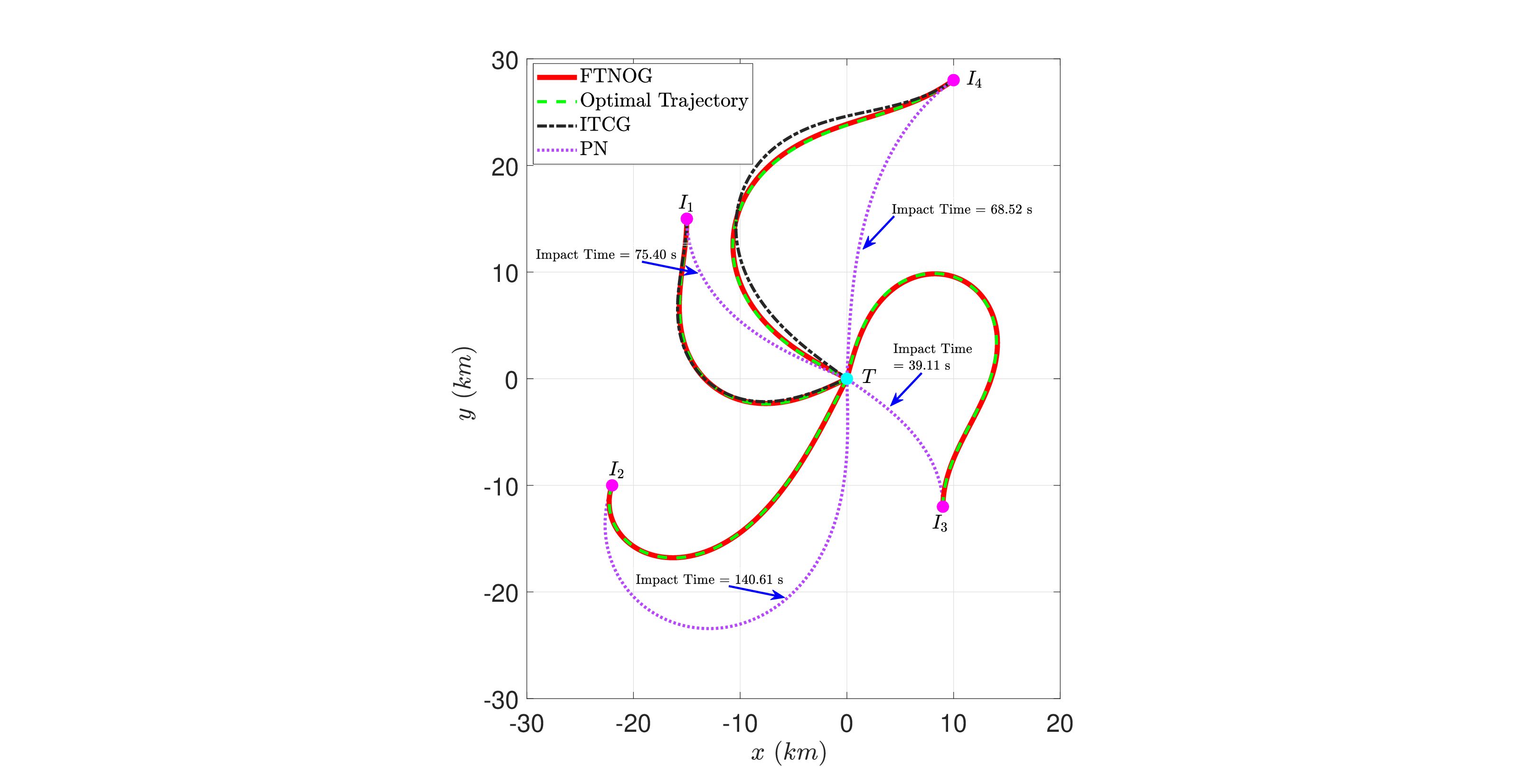}
\caption{Trajectories for cooperative interception of a stationary target.}
\label{Fig:Cooperative_1}
\end{figure}

\begin{figure}[!htp]
\centering
\begin{subfigure}[t]{7cm}
\centering
\includegraphics[width = 8cm]{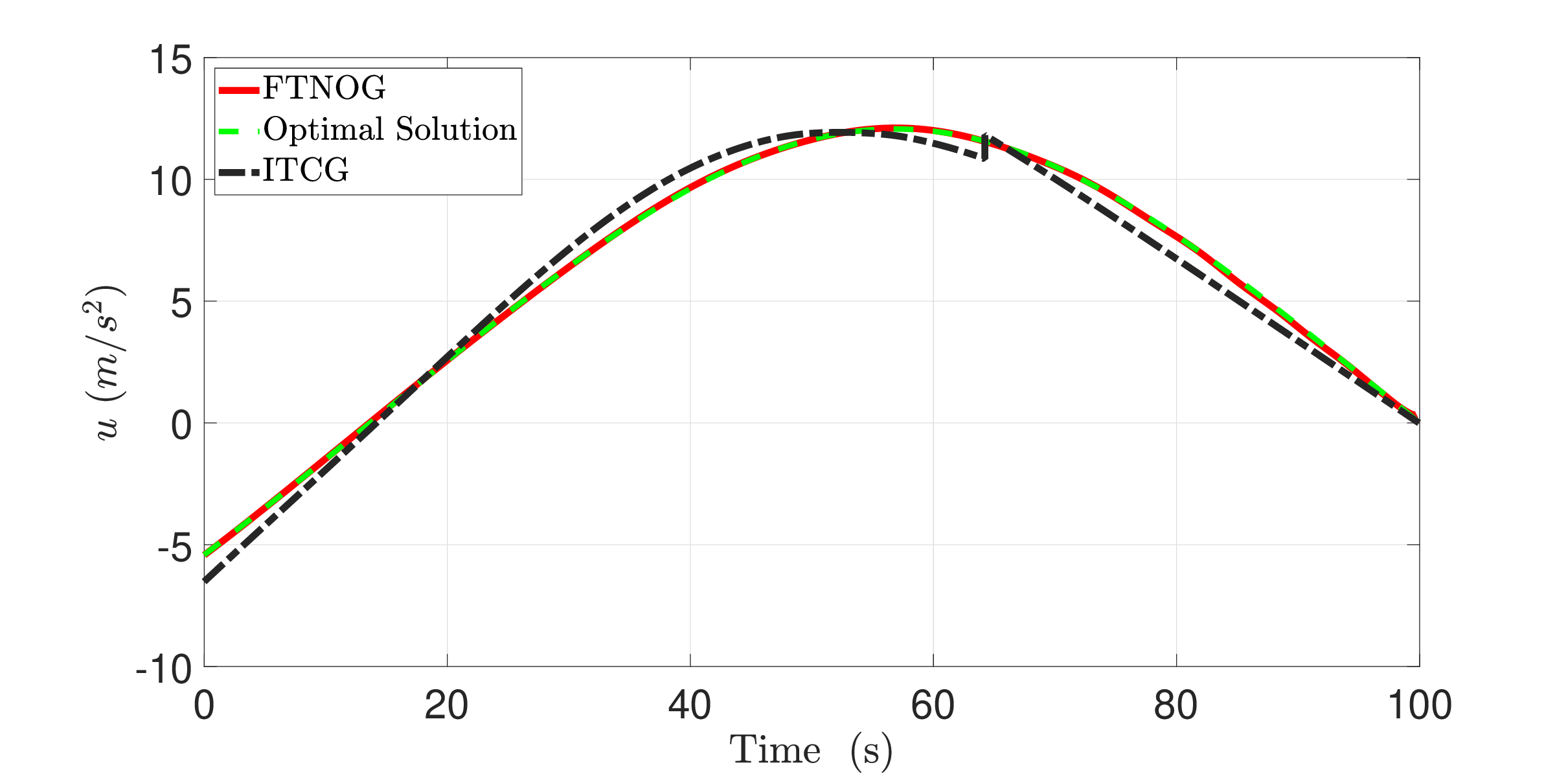}
\caption{Interceptor $\#$1}
\label{Fig:cooperative_control_1}
\end{subfigure}
~~~~~
\begin{subfigure}[t]{7cm}
\centering
\includegraphics[width = 8cm]{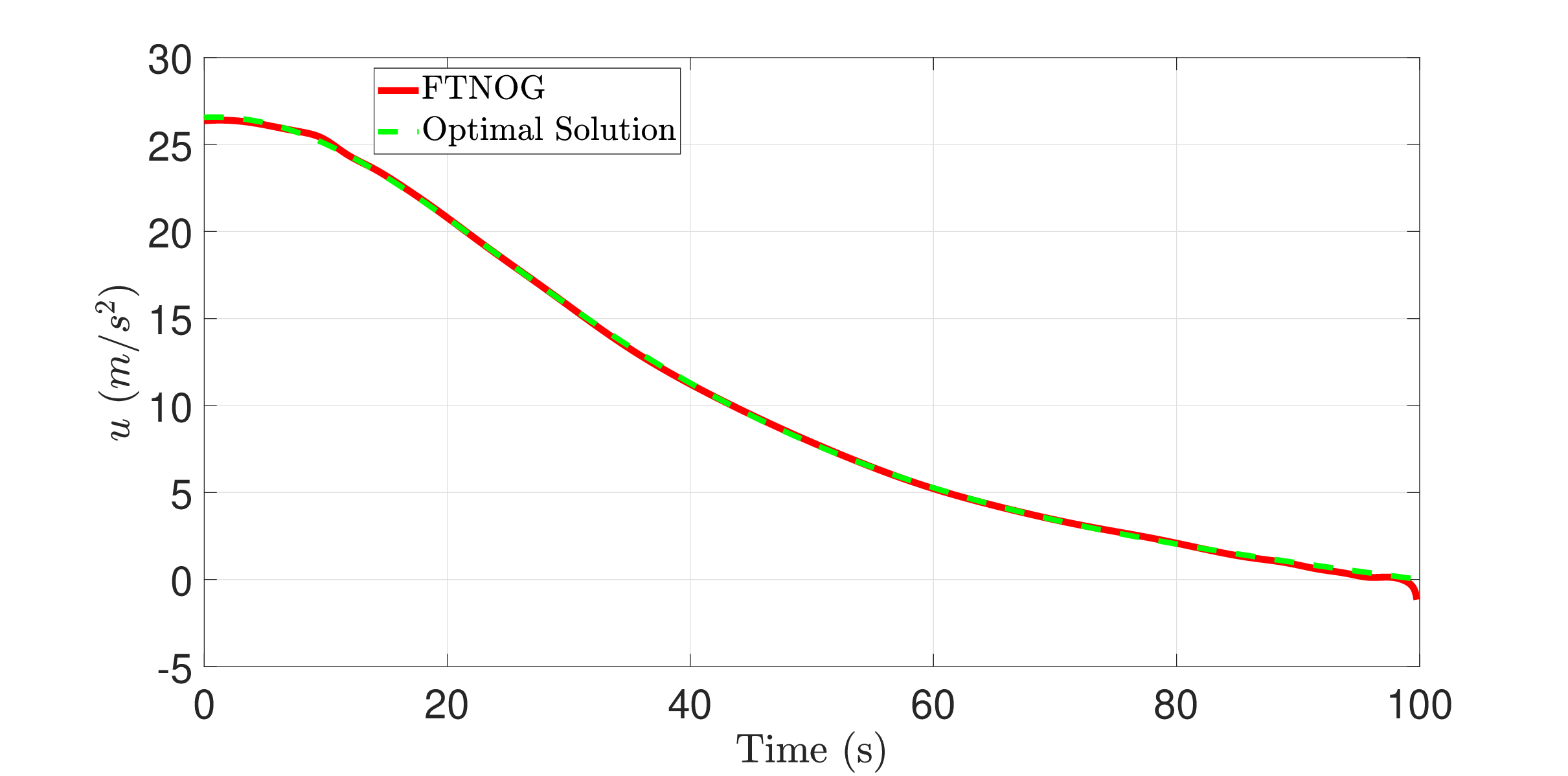}
\caption{Interceptor $\#$2}
\label{Fig:cooperative_control_2}
\end{subfigure}\\
\begin{subfigure}[t]{7cm}
\centering
\includegraphics[width = 8cm]{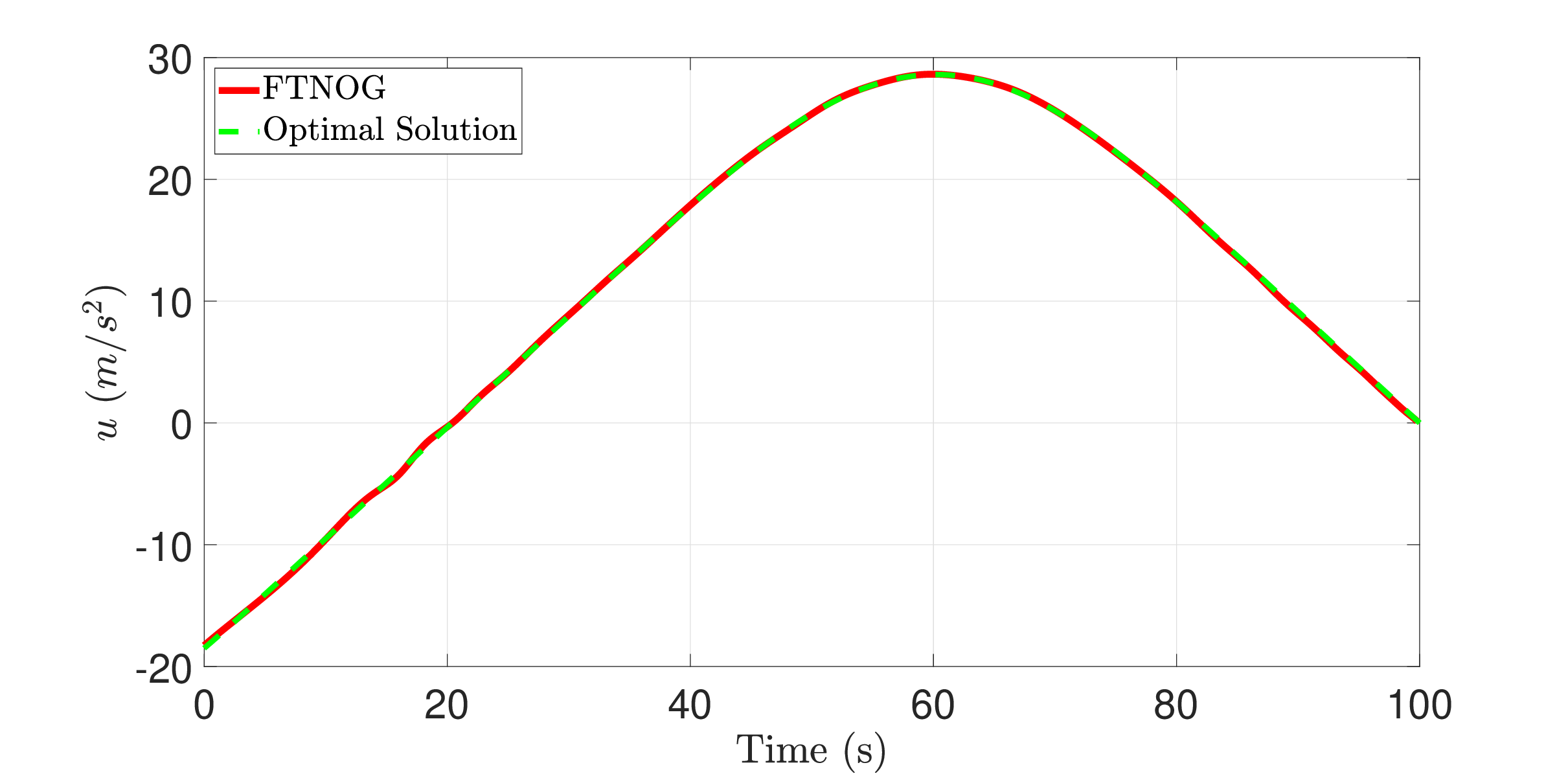}
\caption{Interceptor $\#$3}
\label{Fig:cooperative_control_3}
\end{subfigure}
~~~~~
\begin{subfigure}[t]{7cm}
\centering
\includegraphics[width = 8cm]{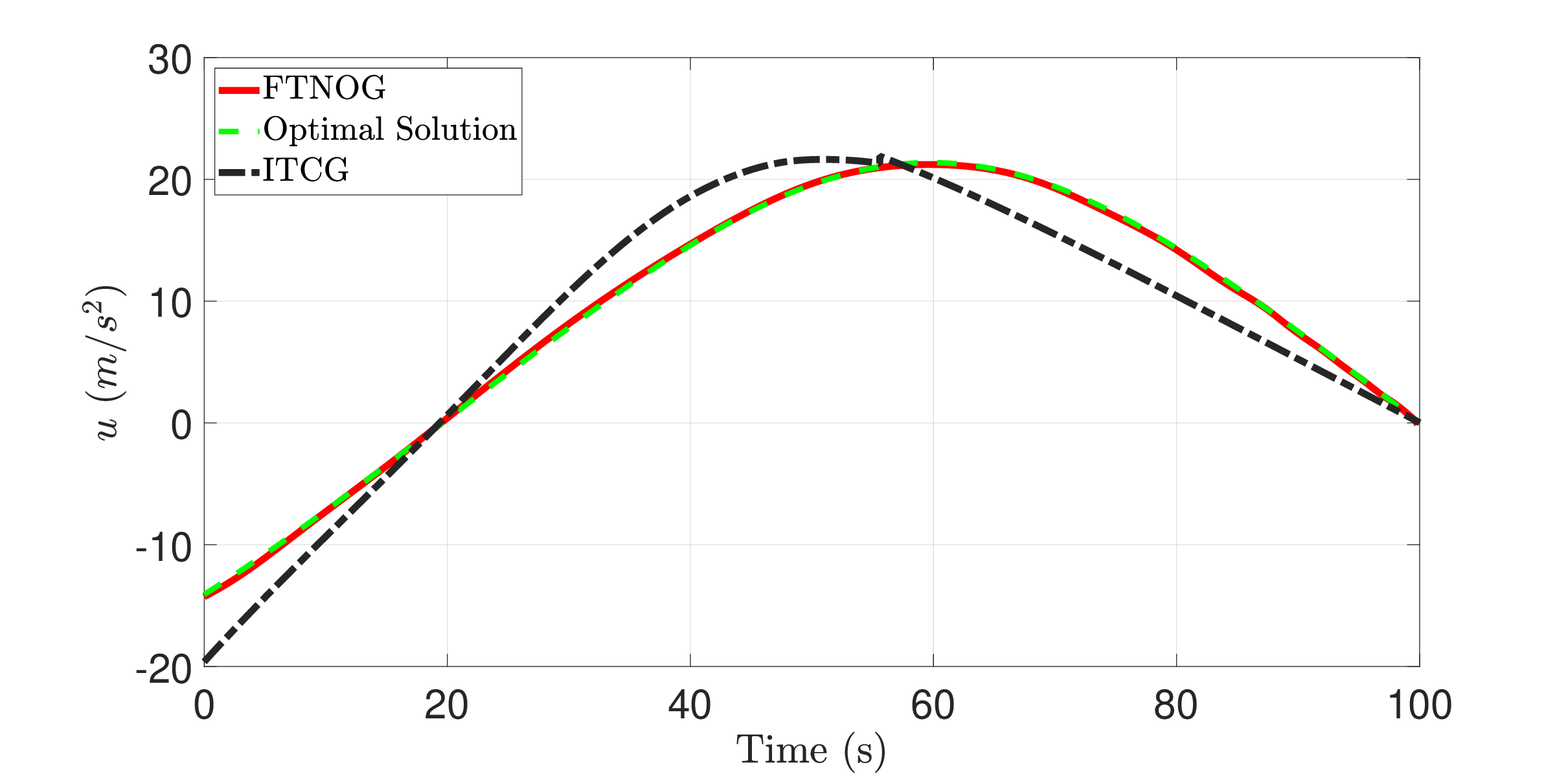}
\caption{Interceptor $\#$4}
\label{Fig:cooperative_control_4}
\end{subfigure}
\caption{Control profiles for trajectories for cooperative interception of a stationary target.}
\label{Fig:cooperative_control_profile}
\end{figure}

In this subsection, an example of salvo attack will be presented to demonstrate the developed FTNOG.  We assume that there are four interceptors. The initial state for each interceptor is given as
\begin{align}
\text{Interceptor \#1}:& \ \ \ (x_0,y_0,\theta_0) = (-15000 ~\text{m}, 15000 ~\text{m},-\frac{\pi}{2} )\nonumber\\
\text{Interceptor \#2}:& \ \ \ (x_0,y_0,\theta_0) = (-22000 ~\text{m},-10000 ~\text{m},-\frac{11\pi}{18} )\nonumber\\
\text{Interceptor \#3}:& \ \ \ (x_0,y_0,\theta_0) = (9000 ~\text{m}, -12000 ~\text{m},\frac{\pi}{2} )\nonumber\\
\text{Interceptor \#4}:& \ \ \ (x_0,y_0,\theta_0) = (10000 ~\text{m}, 28000 ~\text{m}, -\frac{4\pi}{5} )\nonumber
\end{align}
Let $I_i$ denote as interceptor $\# i$, and denote by $V_1$ = $300$ m/s, $V_2$ = $350$ m/s, $V_3$ = $400$ m/s, and  $V_4$ = $450$ m/s the speeds of  $I_1$,   $I_2$, $I_3$, and $I_4$, respectively. The common desired impact time is set as $100$ sec. Then, ITCG, FTNOG, and optimization methods are employed to generate the trajectories, as shown in Fig.~\ref{Fig:Cooperative_1}. The control profiles along trajectories related to different methods are plotted in Fig.~\ref{Fig:cooperative_control_profile}. The ITCG is switched to PN with navigational gain of $3$ once the error of time-to-go estimation reaches zero. In contrary, using the trained network to generate the FTNOG does not require to change to PN in all the engagement.
It is clearly seen from Fig.~\ref{Fig:Cooperative_1} and Fig.~\ref{Fig:cooperative_control_profile}  that the trajectory from FTNOG almost coincides with the optimal trajectory for each interceptor, and the control profiles of FTNOG are close to the controls along optimal trajectories. The discontinuity of control profile indicates that the ITCG is switched to PN.

For Interceptor $\# 1$,  the impact time needed by PN is $75.40$ sec. We can see from  Fig.~\ref{Fig:cooperative_control_1} that the control profile of ITCG is a bit different from the optimal control.
 The control effort related to the FTNOG is closer to the optimal solution than ITCG, as shown by the values in Table.~\ref{Table:cooperative1_control_effort}. It is worth noting that the control effort by PN is smaller than that by optimal control for each interceptor. However, the final impact time for PN cannot be controlled. 
\begin{table}[!htp]
\centering
\caption{The values of control effort used by different guidance laws for salvo attack.}
\begin{tabular}{ccccc}
\hline
      & $J_{ITCG}$~($m^2/s^3$)     & $J_{PN}$~($m^2/s^3$)  & $J_{FTNOG}$~($m^2/s^3$)     & $J^*$~($m^2/s^3$)   \\ 
\hline
Interceptor 1 & $3.1055\times10^3$  &$1.1610\times10^3$  &$3.0978\times10^3$   &$3.0916\times10^3$  \\ 
\hline
Interceptor 2 & NA                  &$6.9592\times10^3$  &$9.4803\times10^3$   &$9.4638\times10^3$  \\
\hline
Interceptor 3 & NA                  &$6.4474\times10^3$  &$1.5838\times10^4$   &$1.5813\times10^4$  \\
\hline
Interceptor 4 & $9.7842\times10^3$  &$1.8474\times10^3$  &$9.4396\times10^3$   &$9.4364\times10^3$  \\
\hline
\label{Table:cooperative1_control_effort}
\end{tabular}
\end{table}

\begin{figure}
\centering
\includegraphics[width = 1\linewidth]{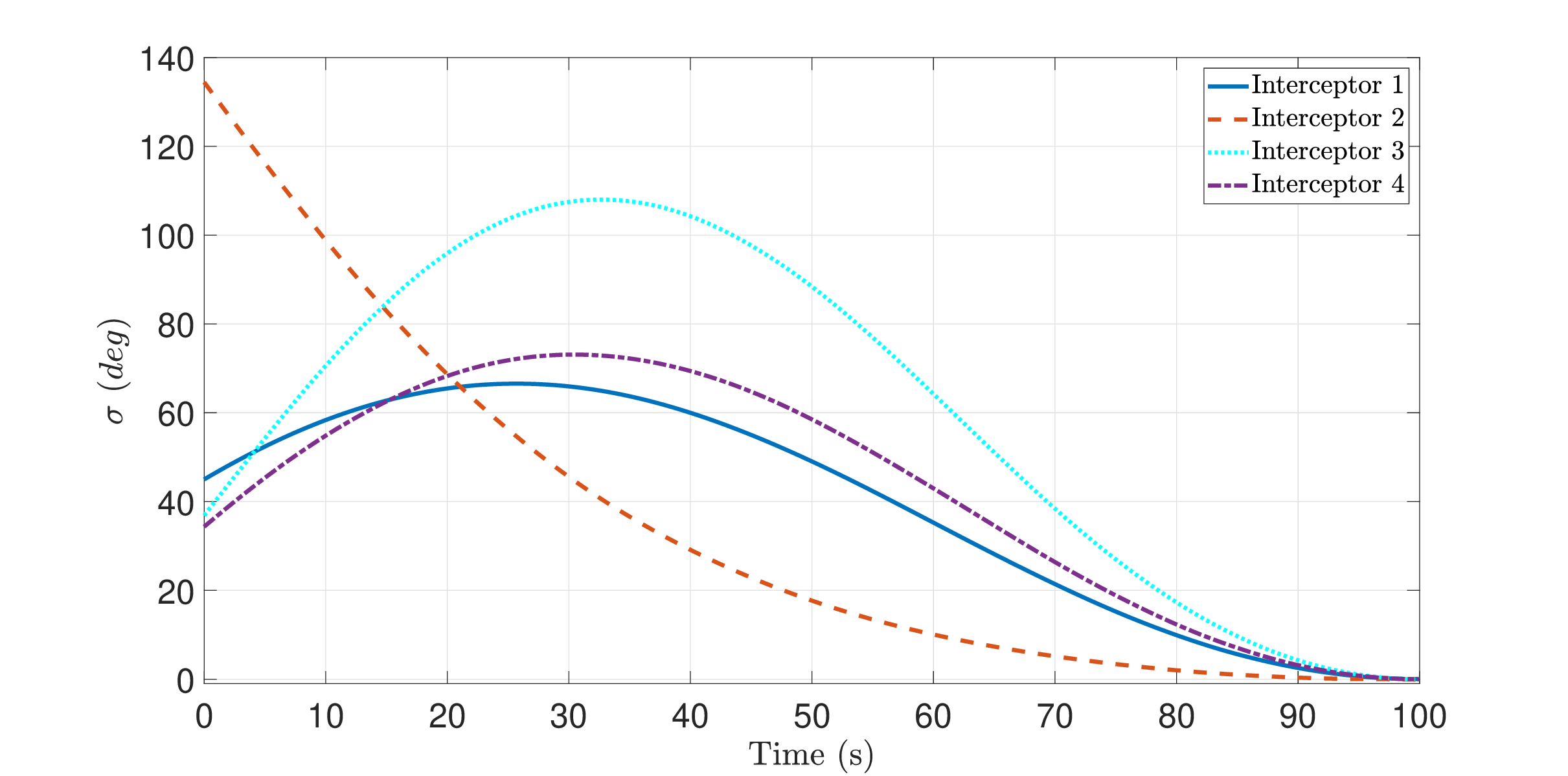}
\caption{Look angle profiles along trajectories related to the FTNOG.}
\label{Fig:Cooperative_sigma}
\end{figure}

Regarding Interceptor $\# 2$, it takes $140.61$ sec to hit the target if  PN is used. According to \cite{Jeon:2006}, the estimated impact time for PN is $217.37$ sec, indicating that the desired impact time is less than the estimated impact time. As a result, the ITCG does not apply to guiding Interceptor $\# 2$ to the target with a duration of $100$ sec. In contrary, the FTNOG is still applicable, providing quite accurate solutions, as shown by control profile in  Fig.~\ref{Fig:cooperative_control_2} and control effort in Table.~\ref{Table:cooperative1_control_effort}.
If using PN to guide Interceptor $\# 3$, it takes $39.11$ sec to impact the target. In this case, the ITCG does not apply either. It is seen from Fig.~\ref{Fig:Cooperative_1} that the FTNOG guides Interceptor $\# 3$ to the target at the desired impact time successfully, by  producing a much longer trajectory. To guide Interceptor $\#4$ to impact the target, the engagement duration required by PN is $68.52$ sec. We can see from Table \ref{Table:cooperative1_control_effort} that the control effort required by FTNOG is close to the optimal solution, and smaller than the ITCG. The profiles of look angles along FTNOG-related trajectories are presented in Fig.~\ref{Fig:Cooperative_sigma}, showing that the optimality condition in Lemma \ref{LE:cos_sigma=1} is met for each trajectory.

\section{Conclusions}\label{SE:conclusions}

Optimal guidance in the nonlinear setting was studied in the paper. Necessary conditions from PMP were analyzed so that a new optimality condition was established. This new optimality condition implies that the look angle along any optimal trajectory can never be zero or $\pi$ (cf. Lemma \ref{LE:cos_sigma=1}), and it also indicates that the optimal control cannot change its sign more than once (cf. Lemma \ref{LE:control_twice}).  By applying the necessary conditions from PMP as well as the new optimality condition,  a parameterized  system was formulated so that its solution space was the same as the solution space of the nonlinear optimal interception problem. In addition, the parameters were bounded according to the symmetric properties of optimal trajectories. All the theoretical developments allowed simply propagating the parameterized system with the parameters in a bounded set to generate enough sampled data for the mapping from current state and time-to-go to the optimal control command (cf. Procedure \ref{algo1}). As a result, a feedforward artificial neural network  could be trained by the sampled data to represent the optimal control command, which eventually was able to generate the FTNOG within $0.1$ millisecond on a usual computation platform. Finally, all the theoretical developments were verified by some numerical simulations, and the performance of the FTNOG was demonstrated by comparing with ITCG and optimal control. 

\section*{Acknowledgement}

This research was supported by the National Natural Science Foundation of China under grant Nos. 61903331 and 62088101.

\begin{appendix}
\section{Proofs for Lemmas in Section \ref{SE:properties}}\label{Appendix:A}

Proof of Lemma \ref{LE:cos_sigma=1}. By contradiction, assume that along the extremal trajectory $(x(\cdot),y(\cdot),\theta(\cdot))$  on $[0,t_f]$, there is a time $\bar{t}\in (0,t_f)$ so that the velocity vector $[\cos \theta(\bar{t}),\sin \theta (\bar{t})]$ is collinear with the LOS, i.e., Eq.~(\ref{EQ:cos_sigma=1}) holds. Let $u(t)$ for $t\in [0,t_f]$ be the corresponding extremal control, and let $A$ be the state at $\bar{t}$, i.e., $A=(x(\bar{t}),y(\bar{t}),\theta(\bar{t}))$. 

\begin{figure}[htbp]
 \centering 
\subcaptionbox{\label{1}}{\includegraphics[trim = 2cm 0cm 2cm 0cm, width = .32\linewidth]{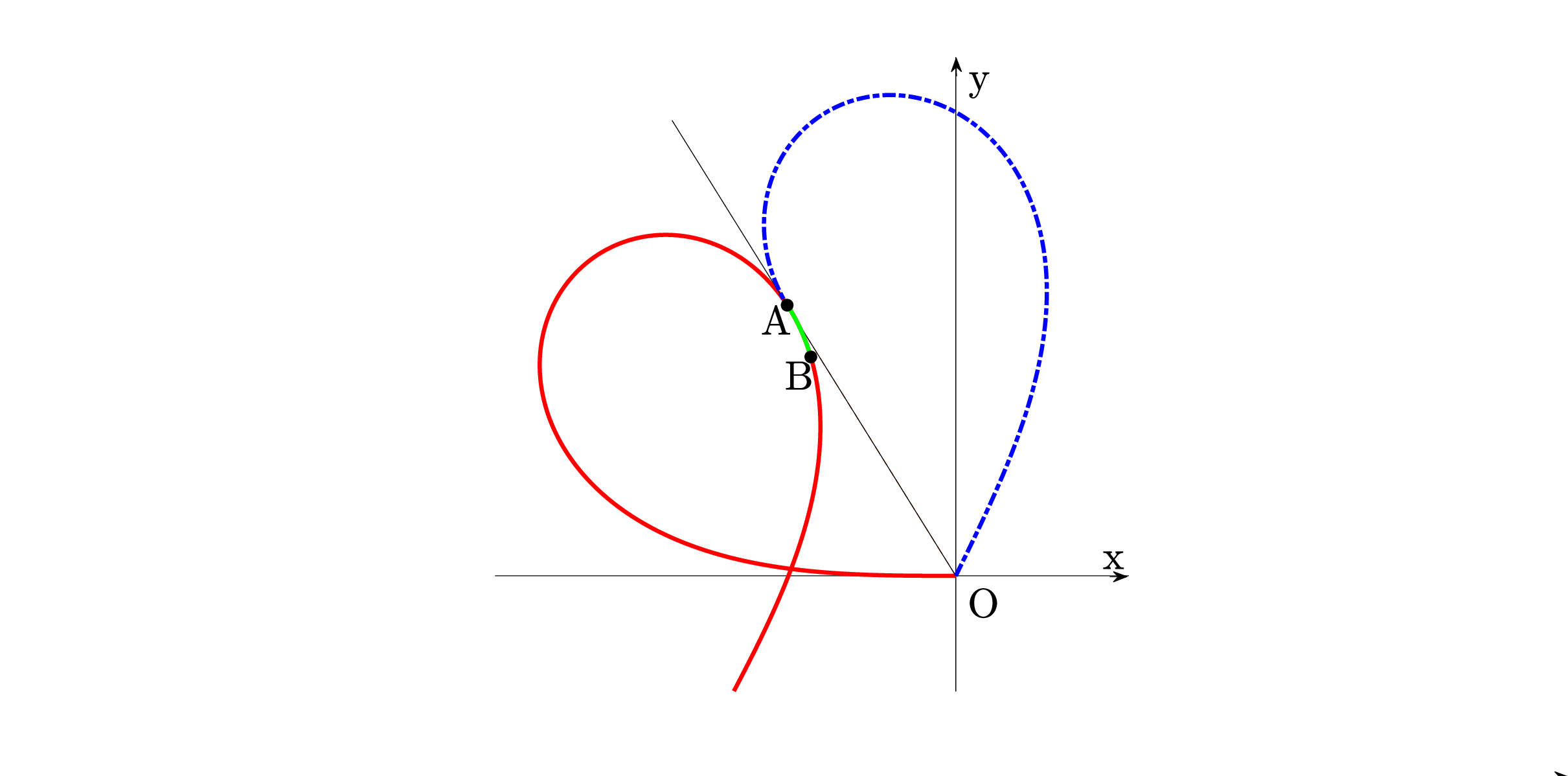}}
\subcaptionbox{\label{2}}{\includegraphics[width = .32\linewidth]{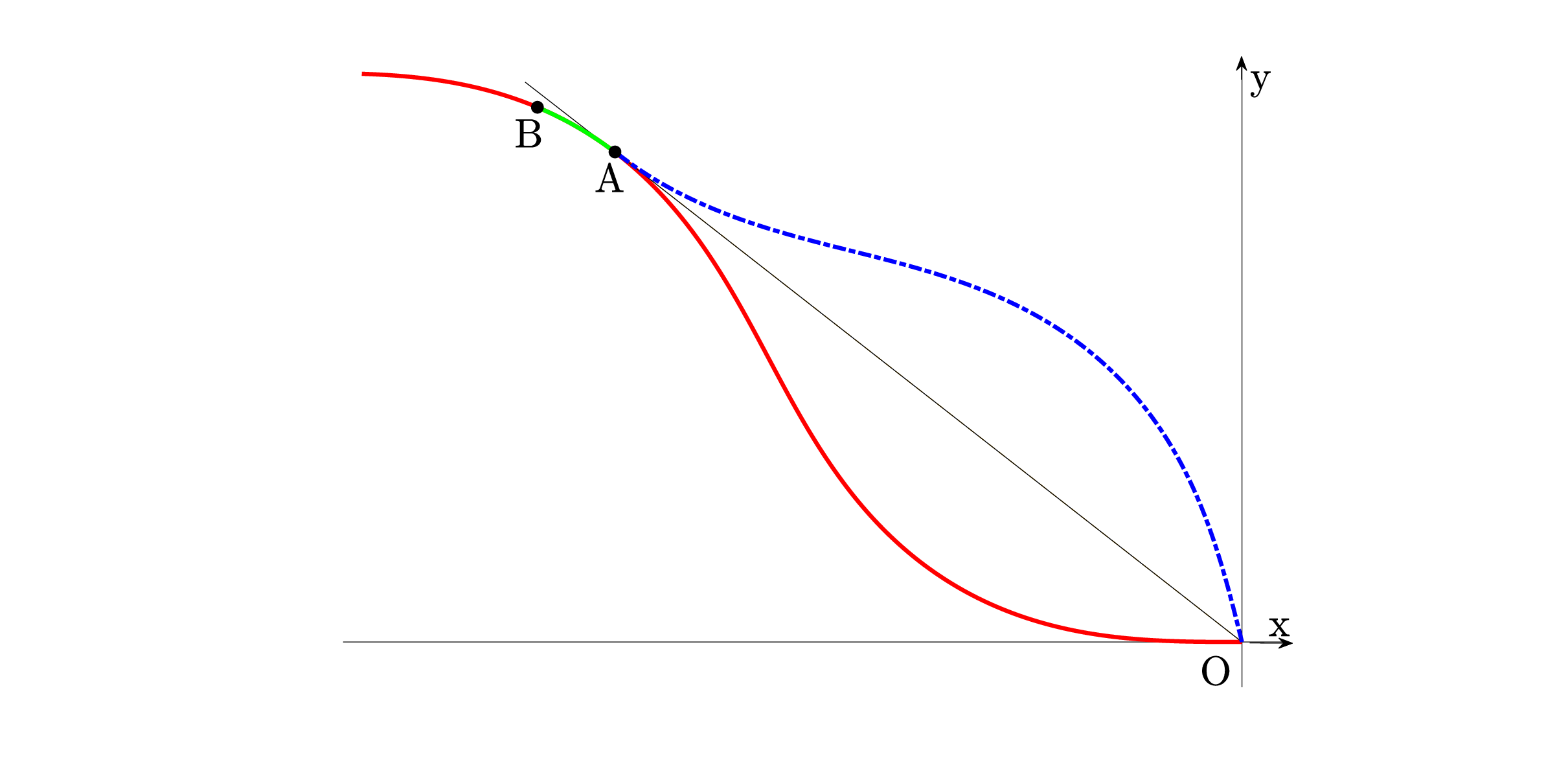}}
\subcaptionbox{\label{3}}{\includegraphics[trim = 3cm 0cm 3cm 0cm, width = .32\linewidth]{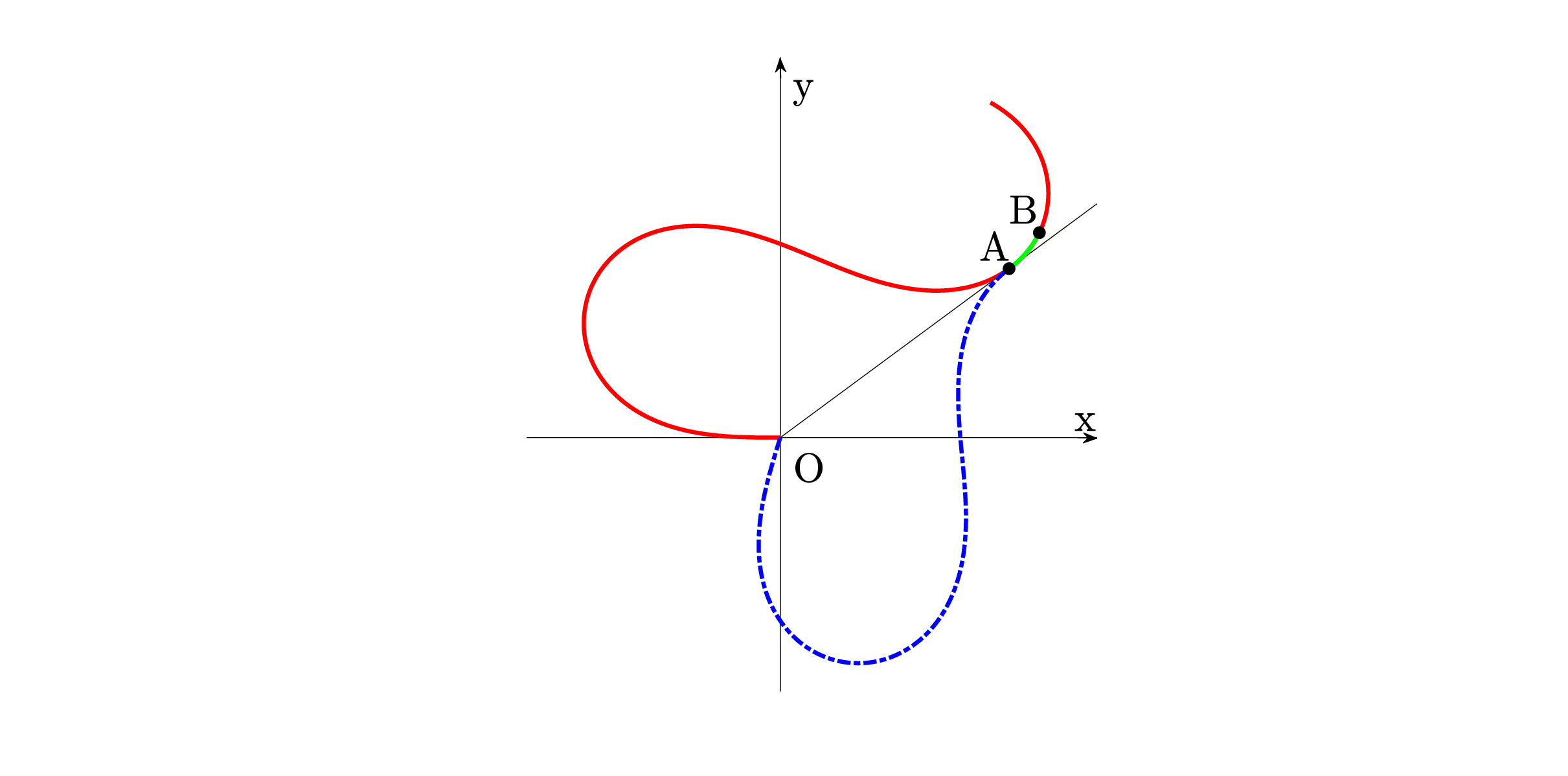}}
\caption{Geometry for existing points with $\sigma = 0$ or $\pi$.}
\label{Fig:optimality}
\end{figure}

Set $\hat{u}(t) = - u(t)$ for $t\in [0,t_f]$. Then, we have that $\hat{u}(t)$ for $t\in [\bar{t},t_f]$ is the extremal control for the interceptor from $A$ to the origin, as shown by the symmetric paths  in Fig.~\ref{Fig:optimality}, where the solid curves denote the extremal trajectory $(x(\cdot),y(\cdot),\theta(\cdot))$ on $[0,t_f]$, and the dashed curves denote the extremal trajectory associated with $\hat{u}(t) = - u(t)$ for $t\in [\bar{t},t_f]$. 

Let us choose a time $\tilde{t}\in [0,\bar{t})$, and let $B$ denote the state at $\tilde{t}$. Denote by $\gamma$ the piece of extremal trajectory $(x(\cdot),y(\cdot),\theta(\cdot))$ from $B$ to the origin, and denote by $\hat{\gamma}$ the smooth concatenation of extremal trajectory $(x(\cdot),y(\cdot),\theta(\cdot))$ from $B$ and $A$ and the extremal trajectory of dashed curve. Then, it is clear that the cost for the interceptor to move along ${\gamma}$ is the same as  the cost of along $\hat{\gamma}$. However, along the new extremal trajectory $\hat{\gamma}$ the control is discontinuous at $A$. This contradicts with the necessary condition in Eq.~(\ref{EQ:u=p_theta}) that the control is continuous. Thus, there is another extremal trajectory from $B$ to the origin with a duration of $t_f - \tilde{t}$ so that the cost is smaller, completing the proof.  
$\square$

Proof of Lemma \ref{LE:control_twice}. 
By contradiction, assume that along an extremal control $u(t)$ for $t\in [0,t_f]$, there are two different times $t_1\in (0,t_f)$ and $t_2\in (0,t_f)$ so that $u(t_1) = u(t_2) = 0$. This indicates $p_{\theta}(t_1) = p_{\theta}(t_2) = 0$. According to Eq.~(\ref{EQ:p_theta}), we have that the three points $(x(t_1),y(t_1))$, $(x(t_2),y(t_2))$, $(0,0)$ lie on the same straight line. Thus, there exists a time $\bar{t}\in (t_1,t_2)$ so that Eq.~(\ref{EQ:cos_sigma=1}) holds. Then, according to Lemma \ref{LE:cos_sigma=1}, the extremal trajectory  $(x(t),y(t),\theta(t))$ for $t\in[0,t_f]$ is not optimal if $u(t_1) = u(t_2) = 0$. This completes the proof.
$\square$

Proof of Lemma \ref{LE:opposite}. Let $(r(t),\sigma(t))$ for $t\in [0,t_g]$ be an optimal solution in the polar frame so that $(r(0),\sigma(0)) = (r_0,\sigma_0)$ and $r(t_g) = 0$, and let $u(t)$ for $t\in [0,t_g]$ be the corresponding optimal control. Then, we have $u(0) = C(r_0,\sigma_0,t_g)$. It is apparent that $(r(t),-\sigma(t))$ for $t\in [0,t_g]$ is an optimal solution in the polar frame, and $-u(t)$ for $t\in [0,t_g]$ is the optimal control associated with $(r(t),-\sigma(t))$, indicating $-u(0)=C(r_0,-\sigma_0,t_g)$. This completes the proof. 
$\square$

\section{Proofs for  Lemmas in Section \ref{SE:Real}}\label{Appendix:B}

Proof of Lemma \ref{LE:h=u}. 
By the definition of $\mathcal{F}_P$, there exists an optimal solution $(r(t),\sigma(t))$, $t\in [0,t_f]$, of the OCP so that
\begin{align}
\begin{split}
r_g = &\ r(t_f - t_g)\\
\sigma_g = &\ \sigma(t_f - t_g)\\
C(r_g,\sigma_g,t_g) =&\  u(t_f - t_g)
\end{split}\label{EQ:LE:h=u1}
\end{align}
where $u(t)$ is the optimal control associated with $(r(t),\sigma(t))$. Up to now, it has been clear that in order to prove this lemma, we just need to prove that there exists $(\alpha,\beta)\in \Lambda$ so that 
\begin{align}
\begin{split}
r(t_f - t) =&\ R(t,\alpha,\beta)\\
\sigma(t_f - t) = &\ \Sigma(t,\alpha,\beta)\\
u(t_f-t) =&\ U(t,\alpha,\beta)
\end{split}
\label{EQ:LE:h=u2}
\end{align}
Notice  that for any optimal trajectory $(r(t),\sigma(t))$ associated with optimal control $u(t)$ for $t\in [0,t_f]$, there exists an optimal trajectory $(x(t),y(t),\theta(t))$ associated with optimal control $u(t)$ for $t\in [0,t_f]$ with $( {x}(t_f), {y}(t_f), {\theta}(t_f)) = (0,0,0)$ so that 
\begin{align} 
r(t) =&\ \sqrt{ {x}(t)^2 + {y}(t)^2}\\
\cos \sigma (t) =&\ \frac{ {x}(t)\cos  {\theta}(t) +  {y}(t) \sin  {\theta}(t)}{\sqrt{ {x}(t)^2 +  {y}(t)^2}}
\end{align}
Then, combining Eq.~(\ref{EQ:sigma_t}), Eq.~(\ref{EQ:U}), and the definition of the parameterized system in Eq.~(\ref{EQ:Sigma1}), we have that there exists $(\alpha,\beta)\in \Lambda$ so that Eq.~(\ref{EQ:LE:h=u2}) holds. This completes the proof. $\square$


Proof of Lemma \ref{LE:r_sigma_in_F}. By definition, we have that for any $(\alpha,\beta)\in \Lambda$, the trajectory $(R(t,\alpha,\beta),\Sigma(t,\alpha,\beta)$ for $t\in [0,T(\alpha,\beta)]$ is an optimal solution, and $U(t,\alpha,\beta)$ is the corresponding optimal control. Thus, we have that the equations in Eq.~(\ref{EQ:LE:r_sigma_in_F}) holds for any $t\in [0,T(\alpha,\beta)]$, completing the proof. $\square$

Proof of Lemma \ref{LE:Norm}. Since $(\bar{r}(t),\bar{\sigma}(t))$ for $t\in [0,\bar{t}_g]$ is an optimal trajectory of an interceptor with speed of $\bar{V}$, we have that the kinematics is
\begin{align}
\begin{cases}
\frac{\text{d}\bar{r}({t})}{\text{d}{t}} = - \bar{V} \cos \bar{\sigma}({t}),\\
\frac{\text{d}\bar{\sigma}({t})}{\text{d}{t}} = \bar{V} \frac{\sin\bar{\sigma}({t})}{\bar{r}({t})} - \frac{\bar{u}({t})}{\bar{V}}\\
\end{cases}
\label{EQ:Sigma_s_bar}
\end{align}
For any $t_g \in (0,\bar{t}_g)$, set $ \tau = t \frac{t_g}{ \bar{t}_g}$, $r(\tau) = \bar{r}(\tau \bar{t}_g/t_g)\frac{t_g}{\bar{t}_g \bar{V}}$, and $\sigma(\tau) = \bar{\sigma}(\tau \bar{t}_g/t_g)$. Then, according to Eq.~(\ref{EQ:Sigma_s_bar}), we have
\begin{align}
\begin{cases}
\frac{\mathrm{d}}{\mathrm{d}\tau} r (\tau) = -\cos \sigma(\tau)\\
\frac{\mathrm{d}}{\mathrm{d}\tau} \sigma(\tau) = \frac{\sin \sigma(\tau)}{r(\tau)} - \bar{u}(t) \frac{\bar{t}_g}{t_g \bar{V}} 
\end{cases}
\end{align}
Since $\bar{u}(t)$ is the optimal control associated with the trajectory $(\bar{r}(t),\bar{\sigma}(t))$ for $t\in [0,\bar{t}_g]$, it follows that $\bar{u}(\tau \bar{t}_g /t_g) \bar{t}_g /(t_g \bar{V})$ is the optimal control associated with the trajectory $(r(\tau),\sigma(\tau))$ for $\tau \in [0,t_g]$. Thus, we have
$C(r(\tau),\sigma(\tau),t_g-\tau) = \bar{u}(t)  \bar{t}_g /(t_g \bar{V})$
which further indicates 
$\bar{u}(t) =C(r(\tau),\sigma(\tau),t_g -\tau) \frac{t_g \bar{V}}{\bar{t}_g } $.
 This completes the proof.
 $\square$
\end{appendix}

\bibliographystyle{unsrt}  
\bibliography{references}

\end{document}